\documentstyle{amsppt}
\magnification = \magstep 1
\pageheight{7 in}
\refstyle{A}
\input xypic
\NoRunningHeads

\define\C{{\Cal C}}

\define\M{{\Cal M}}
\define\di{\operatorname{dim}}
\define\Q{{Q}}
\define\De{{\Cal D}}
\define\CH{{\Cal C}{\Cal M}}
\define\CR{\underline {CR}}
\define\r{\operatorname{real}}

\define\Ho{\operatorname{Hom}}
\define\Pe{{\Cal Perv}}
\define\IC{{\Cal I}{\Cal C}}

\define\En{\operatorname{End}}
\define\Img{\operatorname{Im}}
\define\Ker{\operatorname{Ker}}
\define\Cok{\operatorname{Cok}}
\def\scirc{\scriptstyle \circ \textstyle}

\def\Hom{\operatorname{Hom}}

\topmatter
\title
Motivic decomposition and intersection Chow groups I
\endtitle
\author
Alessio Corti\\Masaki Hanamura
\endauthor
\address D.P.M.M.S., 16 Mill Lane, Cambridge CB2 1SB, UK \endaddress
\email a.corti\@pmms.cam.ac.uk \endemail
\address Dept. of Mathematics, Kyushu University, Hakozaki,
Fukuoka 812--12, Japan \endaddress
\email hanamura\@math.kyushu-u.ac.jp \endemail
\date 13-2-98 \enddate
\abstract For an arbitrary quasiprojective variety $S$, defined over a field
$k$, assumed to be finitely generated over its prime field, we define a 
category $CH \M (S)$ of pure Chow motives over $S$. Assuming conjectures of
Grothendieck and Murre, we prove that the decomposition theorem holds
in $CH \M (S)$. As a consequence, the intersection complex $I S$ of $S$
makes sense as an object of $CH \M (S)$. Part II will give an unconditional
definition of ``intersection'' Chow groups and study some of their properties.
\endabstract
\toc
\specialhead 1. Introduction \endspecialhead
\specialhead 2. Pure motives over a base \endspecialhead
\specialhead 3. Standard conjectures and canonical filtrations \endspecialhead
\specialhead 4. Grothendieck motives over a base, semisimplicity and 
decomposition \endspecialhead
\specialhead 5. Filtrations for quasiprojective varieties \endspecialhead
\specialhead 6. Decomposition in $CH \M (S)$ \endspecialhead
\specialhead References \endspecialhead
\endtoc
\endtopmatter
  
\document

\head 1. Introduction \endhead

We will begin explaining our original motivation for writing this paper,
then move on to summarise our main results and outline the contents of
each section. 

\subhead Motivation \endsubhead

\cite{GoMa} introduced the intersection cohomology groups
$$IH^i\bigl(X({\Bbb C}), {\Bbb Q}\bigr)$$
of a (singular) algebraic variety $X$ defined over the field of complex
numbers (the \'etale version of this theory was constructed in \cite{BBD}
for varieties defined over fields finitely generated over the prime field).
We recall the construction of these groups in the beginning of \S 4.
We list some of their properties:

($1$) There is a factorisation 
$$H^i(X, {\Bbb Q})\to IH^i(X, {\Bbb Q})\to H_{2\di X-i}^{BM}(X, {\Bbb Q})$$
of the Poincar\'e map $H^i(X, {\Bbb Q})\to H_{2\di X-i}^{BM}(X, {\Bbb Q})$.

($2$) There is an intersection product
$$IH^i(X, {\Bbb Q})\times IH^j(X, {\Bbb Q}) \to
H_{2\di X -i-j} (X, {\Bbb Q})$$
which is nondegenerate for proper varieties and intersection cycles of 
complementary dimensions.

($3$) Cohomology acts on intersection cohomology
$$H^i(X, {\Bbb Q})\times IH^j(X, {\Bbb Q})
\to IH^{i+j}(X, {\Bbb Q})$$

This research started out as a program (carried out---to an extent---in 
part II) to define a Chow theoretic analogue
$ICH^r (X, {\Bbb Q})$ of the intersection cohomology groups, satisfying 
corresponding properties, namely:

($1^\prime$) There should be a factorisation 
$$CHC^r(X, {\Bbb Q})\to ICH^r(X, {\Bbb Q})\to CH^r(X, {\Bbb Q})$$
of the natural map $CHC^r(X, {\Bbb Q})\to CH^r(X, {\Bbb Q})$. Here 
$CHC^\bullet$ means ``Chow cohomology'', not the 
operational theory of \cite{FM}, 
which does not have a cycle class map \cite{To}, but the theory developed
in \cite{Ha2,5}.

($2^\prime$) There should be an intersection product
$$ICH^r(X, {\Bbb Q})\times ICH^s(X, {\Bbb Q}) \to CH^{r+s}(X, {\Bbb Q})$$

($3^\prime$) Chow cohomology should act on intersection Chow groups
$$CHC^r(X, {\Bbb Q})\times ICH^s(X, {\Bbb Q})
\to ICH^{r+s}(X, {\Bbb Q})$$

Moreover, there should be a cycle class map $cl:
ICH^r(X, {\Bbb Q}) \to IH^{2r}(X, {\Bbb Q})$ and ($1^\prime -3^\prime$) should
be compatible, via the cycle class maps, with ($1-3$).

This seemed to us an important and interesting question. Since their discovery,
$IH^\bullet \bigl(X({\Bbb C}), {\Bbb Q}\bigr)$, resp. $IH^\bullet 
(X \times \overline k, 
{\Bbb Q}_l)$ have been shown to carry a natural pure Hodge structure,
resp. pure Galois module structure \cite{SaM}, resp. \cite{BBD}. These 
structures are only known to arise in nature as
the cohomology groups of a pure Chow motive (direct summands of cohomology
groups of a smooth algebraic variety. We summarise Grothendieck's definition
of pure motives in the beginning of \S 2). If this ``intersection motive'' can 
be identified and constructed, one can then take its Chow groups.

Where to look for such a motive? Let us assume for simplicity that $X$ is
a variety with a single isolated singularity, having a resolution
$f:Z \to X$ introducing a single smooth exceptional 
divisor $E$. According to the decomposition theorem \cite{BBD, SaM} 
$$Rf_\ast {\Bbb Q}_Z={\Cal IC}_X \oplus V$$
where ${\Cal IC}_X$---in Borel's notation \cite{Bo}---is the intersection 
complex of $X$ (its hypercohomology  
groups are the $IH^i(X, {\Bbb Q})$) and $V$ is a sheaf supported on the 
singular point. For the reader's convenience, we will recall the definition of
intersection complexes and the statement of the decomposition theorem early 
on in \S 4. 
If the decomposition theorem is to hold for motives, we
expect the ``intersection motive'' $IX$ of $X$ to be a direct summand
of the motive $hZ$ of the desingularisation. If the decomposition theorem
is to hold for motivic sheaveas on $X$, we expect $IX$ to be of the 
form $(Z, i_\ast P)$ where $P \in CH_{\di X} (E \times E)$ is a projector
in the correspondence ring of $E$ and $i:E \times E \to Z \times Z$ the
inclusion. It is relatively easy to figure out 
what the cohomology class of $P$ must be: $E$ is naturally polarised 
by the dual of its normal bundle, and $P$ must induce the Hodge $\Lambda$
operator (see the beginning of \S 3 if you don't remember this)
relative to this polarisation. In other words, in this 
particular case, {\it the 
motivic decomposition theorem is equivalent to the standard
conjecture of Lefschetz type for $E$} (actually, there is a further quite
subtle problem which will not be discussed here in the introduction, to 
justify that the {\it Chow} motive $(Z,i_\ast P)$ is independent on the 
choice of $P$). On the 
one hand this is quite disappointing: there seems no way to have a reasonable
theory without proving the standard conjecture. On the other hand, we 
are at least able now to place the original question in its proper framework.

\subhead Main results \endsubhead

This paper is not actually concerned with intersection Chow groups, we
plan to do those in part II. Our first result is

\proclaim{Theorem 1} (See \S 2 for
precise statements) let $S$ be an arbitrary quasiprojective
variety, defined over a field $k$, assumed to be finitely generated over its
prime field. There is a category $CH \M (S)$ of pure Chow motives over $S$, 
with realisation in $\De^b_{cc} (S)$, which
is the relative analogue of the category $CH \M$ defined by Grothendieck.
As $CH \M$, $CH \M (S)$ arises from a correspondence category $CH\C (S)$ whose
objects are smooth varieties $X$, together with a projective morphism
$X \to S$, and morphisms are defined as
$$\Ho_{CH \C S} (X,Y)= \oplus CH_{\di Y_\alpha}(X \times_S Y_\alpha)$$
the sum being taken over all irreducible components $Y_\alpha$ of $Y$.
The composition of morphisms uses Fulton and Mac Pherson's refined Gysin
maps.\endproclaim

The construction of $CH \M (S)$ is very easy and generalises an earlier
idea by \cite{DM}. We are convinced that $CH \M (S)$ is a useful language,
see for instance \cite{DM}, \cite{Sc1}.

Our second main result is

\proclaim{Theorem 2} (See \S 6 for the precise statement) Consider, as 
before, a quasiprojective variety $S$ over $k$, finitely generated over
its prime field. Assume that resolutions of singularities exist for varieties
over $k$. Then, assuming Grothendieck's
standard conjectures and Murre's conjecture, a decomposition theorem
holds in $CH \M (S)$ which realises to the (topological)
decomposition theorem in $\De^b_{cc}(S)$ of \cite{BBD}.\endproclaim

The assumption on resolutions of singularities is probably unnecessary:
the modifications of \cite{DJ1--2} are possibly sufficient for our purposes.
We haven't pursued this, since it is hardly the point of the paper.

The result is part of a larger program, due to M. H., to construct
a triangulated category $\De (S)$ of ``mixed motivic sheaves'' on $S$, and
show that, assuming the standard conjectures and Murre's conjecture
(and, in addition, the vanishing conjecture for $K$-groups of Soul\'e
and Beilinson), it possesses the expected $t$-structure.

The decomposition ``theorem'' allows to make sense of intersection complexes
and intersection Chow groups. In part II we shall use the ideas and
result of part I to propose an unconditional definition of intersection
Chow groups, and study some of their properties, sometimes with the aid 
of the conjectures.

We are sure that the reader will perceive our liberal use of various 
long standing conjectures to be a significant weakness of our study.
As a partial answer to this possible objection, we would like to
make 2 remarks. First, we are making the point in this
paper, following an insight of M. H., that the standard conjectures, 
which were designed primarily to deal with motives over the point, are
indeed enough to determine the first order (i.e., pure) behaviour
of {\it motivic sheaves}. Second, there are interesting concrete 
contexts where the conjectural assumptions are satisfied, or could 
conceivably be shown. These include families of curves, surfaces, abelian
varieties, toric morphisms. We believe that our theory will prove to be
useful in these situations.

Before giving a quick detailed description of the contents of each section, 
we would like to say that we invested an inordinate amount of time to make 
this paper as self contained as possible. We assume a basic knowledge
of algebraic geometry and intersection theory, as can be accessed through
the 1st 6 chapters of \cite{Fu}, and a working knowledge of Verdier duality
as can be obtained e.g. by looking into \cite{KS}. We shall recall or
summarise everything else we need, and that includes the bivariant theory
of \cite{FM}, the standard conjectures, Murre's conjecture on the 
natural filtration on the Chow groups
of smooth and projective varieties, perverse sheaves, the decomposition 
theorem,
etc. We especially hope that this paper will be accessible to nonexperts
and graduate students seeking an introduction to the field.
 
\subhead Summary of contents \endsubhead

\S 2 is devoted to the construction of $CH \M (S)$ and the study of its first
properties, especially the realisation in $\De^b_{cc} (S)$. In \S 3
we recall the standard conjectures of Grothendieck, Murre's conjecture
on the canonical filtration of the Chow groups of smooth and projective
varieties, and S. Saito's proposed unconditional
definition of this filtration. In \S 4 we define a category $\M (S)$ 
of ``Grothendieck'' motives over $S$ and, assuming the standard conjectures,
we show that it is abelian and semisimple. This is an intermediate
step in the direction of the decomposition theorem in $CH \M(S)$, 
which we finally
prove in \S 6, after extending Saito's filtration to Chow groups
of quasiprojective varieties in \S 5. For more information on the material
covered in the various sections and the logical dependencies between them, 
the reader is invited to consult the short summary that we provide at the 
beginning of each.

\subhead Acknowledgements \endsubhead

The authors would like to thank the mathematical institute of the University
of Warwick for making it possible for them to meet in the Fall 1995 
during the WAG special year. A. C. would like to 
express his
immense gratitude to M. H. for teaching him the vast body of knowledge
comprised in \cite{Ha1--6} and for his infinite patience during the long time
of preparation of the manuscript.

\head 2. Pure motives over a base \endhead

The main object of this section is to generalise Grothendieck's 
construction of
pure motives to the relative situation over an arbitrary quasiprojective
veriety $S$. In short, we will define a category $CH \M (S)$ of Chow motives
over $S$; a way to think of it is the category of pure motivic sheaves over 
$S$. After reviewing our conventions for cohomology theories, we recall
Grothendieck's construction of the category $CH \M$ of Chow motives, then
move on to define $CH \M(S)$ along very similar lines: the new element is
the definition of composition of morphisms in $CH\M (S)$ using Fulton and
Mac Pherson's refined Gysin maps. We close the section with the construction
of a natural realisation functor $CH \M (S) \to \De^b_{cc} (S)$. This
construction is quite technical and we will complete it after a short
summary of the topological bivariant theory of \cite{FM}, 
which is an essential ingredient in the proof: we advise the reader to skip
it on first reading.

\subhead {} Cohomology theories \endsubhead

In this subsection, we explain our notation and 
conventions for cohomology theories.

We will need some properties which are not shared by all Weil cohomology
theories (in particular De Rham theory), especially the bivariant
formalism, but later on, as we progress in the study, we will also need
perverse sheaves and reasonable specialisation properties. For this
reason we will work either with Betti cohomology
or with \'etale cohomology. 

(1) We fix a field $k$, finitely generated over the prime field, and consider 
quasiprojective varieties defined over $k$.

The notation $H^i X$ means either:

$H^i\bigl( X({\Bbb C}), {\Bbb Q} \bigr)$, if $k$ has characteristic 0, where
we always assume to have chosen an embedding $\sigma : k \hookrightarrow
{\Bbb C}$, or 

$H^i(\overline X , {\Bbb Q}_l)$, if 
$\operatorname{char} k \not = l$, where $\overline X = X \otimes
\overline k$.

These are vecor spaces over $Q={\Bbb Q}$ or ${\Bbb Q}_l$, depending 
on the context. The mixed
Hodge structure or mixed Galois module structure on these spaces will
be unimportant for us and, for this reason, we do not keep track of
Tate twists in our notation for cohomology groups.
We denote $H^{BM}_i X$ the Borel-Moore homology theory companion to $H^i X$.
Homology and cohomology are part of the more general
bivariant formalism of \cite{FM}, which we shall summarise below when needed. 
There is a natural cycle class map $cl:CH_rX \to H^{BM}_{2r}X$.

(2) If $S$ is a quasiprojective variety defined over $k$, 
$D^b_{cc} S$ denotes either \linebreak
$D^b_{cc}\bigl(S({\Bbb C}), {\Bbb Q}\bigr)$, the
derived category of cohomologically constructible (for the euclidian topology)
sheaves, or $D^b_{cc}\bigl(\overline S, {\Bbb Q}_l\bigr)$, the category 
constructed in \cite{BBD}. 
This has the 6 operations of Grothendieck,
Verdier duality etc. $Q_S$, resp. $D_S$ will denote the constant sheaf,
resp. dualising sheaf, i.e., $Q_S$ is either ${\Bbb Q}_{S({\Bbb C})}$ 
or ${\Bbb Q}_{\overline S, l}$.
Cohomology, Borel-Moore homology and the
bivariant formalism alluded to above arise
from $D^b_{cc}$ and the 6 operations in a familiar way \cite{FM}, which 
is also briefly recalled below.

(3) In \S 3 we will use the following specialisation property of $H^i X$.
Let $T$ be the spectrum of a discrete valuation ring with residue field
$k$ and quotient field $K$ (both finitely generated 
over their prime field), $0, \eta \in T$ the central and generic point, 
$X \to T$ a smooth and proper
morphism. There is then a specialisation isomorphism
$$hsp: H^iX_\eta \overset \cong \to \to H^i X_0$$
compatible via the cycle class with the specialisation homomorphism for
Chow groups 
$$\diagram
CH^r X_\eta \rto^{csp} \dto_{cl}& CH^r X_0\dto^{cl}\\
H^{2r}X_\eta \rto_{hsp}         & H^{2r}X_0 \enddiagram
$$

\subhead {} Reminder of Grothendieck motives \endsubhead

In this subsection, we give a quick reminder of Grothendieck's 
classical construction of motives (over the point), while also
fixing our notation. This construction is 
in 3 steps: first the construction of a correspondence category,
followed by pseudoabelianisation and the introduction of Tate objects and
twists by them.
 
The standard references for this material are \cite{De}, \cite{Ma}, 
\cite{Sc}.

\subsubhead Correspondences \endsubsubhead

We fix a field $k$, finitely generated over its prime field. 
We consider smooth and projective varieties $X$ over $k$, and denote
$C^i X$ the group of algebraic cycles of codimension $i$ 
on $X$ modulo a suitable equivalence
relation. The examples are:

(1) $C^i X = CH^i X$, the Chow group of cycles modulo {\it rational} 
equivalence. 

(2) $C^i X =  A^i X = H^{2i} X_{alg} = \Img (CH^i X \to H^{2i} X)$ is the 
group of cycles modulo {\it homological} equivalence.

(3) $C^i X =$ cycles on $X$ modulo {\it numerical} equivalence.

We will now construct the categories $C {\Cal C}$ of $C$-{\it correspondences}.

\definition{2.1 Definition}
An {\it object} of $C{\Cal C}$ is a {\it smooth and projective}, not 
necessarily
connected, variety $X$.

{\it Morphisms} in $C {\Cal C}$ are correspondences: 
$$\Ho_{C \C}(X, Y)=\oplus C^{\di X_\alpha}X_\alpha \times Y$$
where $X=\coprod X_\alpha$ is the decomposition of $X$ into its
connected components $X_\alpha$.

Let $u: X_1 \to X_2$ and $v: X_2 \to X_3$ be 
correspondences, let $p_{ij}: X_1\times X_2 \times X_3
\to X_i \times X_j$ be the projection. The {\it composition}
is defined as follows:
$$v \circ u= p_{13 \ast} (p_{23}^\ast v \cdot p_{12}^\ast u)$$
\enddefinition

It is easy to see that, with the above definitions, $C \C$ is an
additive category, with the disjoint union of
varieties being the categorical direct sum.

Since the intersection product for Chow groups is compatible with the 
cup product for cohomology classes, we 
have a forgetful functor $CH \C \to A\C$ from the category of Chow
correspondences to the category of homological correspondences.

\subsubhead Chow and Gro\-then\-dieck motives \endsubsubhead

We first recall the construction of the pseudoabelianisation of 
an additive category. Let ${\Cal A}$ be an additive category,
$A$ an object of ${\Cal A}$. A projector is an arrow $P: A \to A$ such that
$P^2=P$. It is possible to give a categorical definition of the image of
a projector:

\definition{2.2 Definition} The {\it image} of a projector $P: A \to A$
is an object $\Img P$ of ${\Cal A}$, together with a factorisation 
of $P : A \to A$ (commutative diagram)
$$\diagram
A \rrto^P \drto & &A \\
                &\Img P \urto &\enddiagram$$
satisfying the 2 identities
$$\Ho (-, \Img P)=P\circ \Ho(-, A)$$
$$\Ho (\Img P, -)=\Ho(A, -)\circ P$$
${\Cal A}$ is {\it pseudoabelian} if every projector has an image.
\enddefinition

\definition{2.3 Remark} The definition of image of $P$ simply means that, for
all objects $X$, $\Ho (X, I)$ is the image, in the category of {\it abelian
groups} of $P\circ \_$ by means of the following diagram
$$
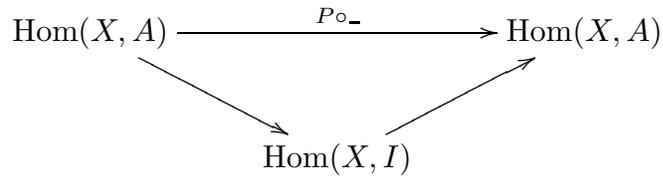
\diagram
\Ho (X,A) \rrto^{P \circ \_ } \drto & &\Ho (X,A) \\
                &\Ho(X,I) \urto &\enddiagram$$
and, at the same time, $\Ho (I, X)$ is the image of $\_\circ P$ by means
of the following diagram
$$
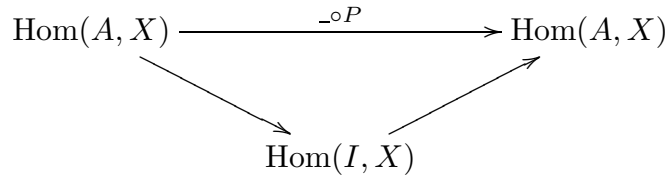
\diagram
\Ho (A,X) \rrto^{\_\circ P} \drto & &\Ho (A,X) \\
                &\Ho (I, X) \urto &\enddiagram$$
\enddefinition

\definition{2.4 Definition} The {\it pseudoabelianisation} of $\Cal A$ is
the category $\tilde {\Cal A}$ defined as follows.

{\it Objects} of $\tilde {\Cal A}$ are pairs $(A, P)$ of an object $A$
of ${\Cal A}$ and a projector $P: A \to A$.

{\it Morphisms} in $\tilde {\Cal A}$ are defined as follows
$$\Ho_{\tilde {\Cal A}} \bigl( (A, P), (B, Q) \bigr)=
Q\circ  \Ho_{\Cal A} (A, B) \circ P$$  
It is a simple observation that this is the same as morphisms $f: A \to B$ in 
${\Cal A}$ such that $f= Q \circ f\circ P$.
\enddefinition

The following result is a formal exercise:

\proclaim{2.5 Theorem} The category $\tilde {\Cal A}$ is pseudoabelian.
There is a natural functor $F: {\Cal A} \to \tilde {\Cal A}$.
Let ${\Cal B}$ be a pseudoabelian category and $G: {\Cal A} \to {\Cal B}$
a functor. Then there exist a unique functor 
$H: \tilde {\Cal A}\to {\Cal B}$ such that $G= H \circ F$.
\endproclaim 

We now define the category $C \M$ of $C$-motives. This is made by taking
the pseudoabelianisation of $C \C$ and then inserting Tate objects and twists
by them:

\definition{2.6 Definition} An {\it object} of $C \M$ is a triple
$$(X, P, r)$$
also denoted $(X,P)(r)$, where $X$ is a 
{\it smooth projective}, not necessarily
connected variety, $P\in \En_{C \C}(X,X)$ a projector, and $r \in {\Bbb Z}$ is
an integer.

{\it Morphisms} in $C \M$ are defined as
$$\Ho_{C \M} \bigl( (X,P,r), (Y,Q,s)\bigr)=
Q \circ \bigl( \oplus C^{\di X_\alpha +s-r} (X_\alpha \times Y)\bigr) \circ P$$
where $X = \coprod X_\alpha$ is the decomposition of $X$ into its connected
components $X_\alpha$. Composition is by means of the same formula 
used for composing correspondences.
\enddefinition

\definition{2.7 Remarks, terminology, notation, etc.}

(1) For $C = CH$, the category $CH \M$ is called the category of {\it Chow}
motives. If $C$ is cycles modulo numerical equivalence, the corresponding 
category of motives is denoted simply $\M$ and called the category of
{\it Grothendieck} motives. As noted below in 3.5, one of the consequences
of the standard conjectures is that $A\M=\M$.

(2) Denoting $\Cal V$ the category of smooth and projective varieties over
$k$, there are natural contravariant {\it cohomological} $h: {\Cal V} \to
CH \M$ and covariant {\it homological} $h^\vee : {\Cal V} \to CH \M$ functors.
As an object, $hX =X$ regarded as a Chow motive, and for a morphism
$f: X \to Y$, $h(f)=cl \Gamma_f^t$ is the cycle class of the transpose
$\Gamma_f^t \subset Y \times X$ of the graph $\Gamma_f$ of $f$.
Similarly, $h^\vee X = \oplus X_\alpha(\di X_\alpha)$, where $X=\coprod
X_\alpha$ is the decomposition in connected components, and
$h^\vee (f)=cl \Gamma_f$.

(3) Similarly,   
there are natural contravariant cohomological $\underline H: 
{\Cal V} \to \M$ and covariant homological 
$\underline H^\vee: {\Cal V} \to \M$ functors.
These are defined in a similar way to $h$ and $h^\vee$.

(4) If $Vec_Q$ is the category of vector spaces over $Q$, where the
cohomology theory $H^i X$ is taking values, there are also
{\it realisation} functors $H^\ast : \M \to Vec_Q$ sending a motive to its
cohomology.
\enddefinition

\subhead {} Pure motives over a base \endsubhead

We extend Grothendieck's construction to the case of varieties over an 
arbitrary
quasiprojective base variety $S$. In doing so, we generalise \cite{DM}.
Realisation functors need more work and we treat them in the next subsection. 

\subsubhead Correspondences over $S$ \endsubsubhead

We fix a field $k$, finitely generated over its prime field. 
We consider quasiprojective varieties $Z$ over $k$, and denote
$C_i Z$ the group of $i$-dimensional algebraic cycles on $Z$ modulo a 
suitable equivalence relation. The examples are:

(1) $C_i Z = CH_i Z$, the Chow group of cycles modulo {\it rational} 
equivalence. 

(2) $C_i Z =  A_i Z = H^{BM}_{2i} X_{alg} = \Img (CH_i X \to H^{BM}_{2i} X)$ 
is the group of cycles modulo {\it homological} equivalence.

Let us fix an arbitrary quasiprojective variety $S$.
We will now construct the categories $C {\Cal C}(S)$ of 
$C$-correspondences {\it over S}.

\definition{2.8 Definition}
An {\it object} of $C{\Cal C}(S)$ is a {\it smooth}, not necessarily
connected, variety $X$, together with a {\it projective} morphism
$f: X \to S$.

{\it Morphisms} in $C {\Cal C}(S)$ are correspondences: 
$$\Ho_{C \C(S)}(X, Y)=\oplus C_{\di Y_\alpha}(X \times_S Y_\alpha)$$
where $Y = \coprod Y_\alpha$ is the decomposition of $Y$ into its
connected components $Y_\alpha$.

The {\it composition} of morphisms is realized with the help of the following
fibre square diagram
$$\diagram 
X\times_S Y \times_S Z \rto \dto & (Y\times_S Z)\times (X\times_S Y) \dto
\\ Y \rto_\delta & Y \times Y \enddiagram$$
For $u: X\to Y$, $v:Y \to Z$ we define the composition
$v \bullet u :X \to Z$ by
$$v \bullet u= p_{XZ \ast} \delta^! (v \times u)$$
where $p_{XZ}$ is the projection on the first and third factor
$$p_{XZ}:X\times_S Y \times_S Z \to X\times_S Z$$
and  $\delta^!$ is Fulton's
refined Gysin map for local complete intersection (lci) morphisms 
\cite{Fu, BFM}. 
We are assuming that $Y$ is
smooth, therefore the diagonal embedding $Y \to Y \times Y$ is lci.
\enddefinition

It is easy to see that, with the above definitions, $C \C(S)$ is an
additive category, with the disjoint union of
varieties being the categorical direct sum. 

Since the intersection product for Chow groups is compatible with the 
cup product for Borel-Moore homology classes \cite{BFM}, we 
have a forgetful functor $CH \C (S) \to A\C (S)$.

\subsubhead Chow and homological motives over $S$ \endsubsubhead

We now define the category $C \M (S)$ of pure $C$-motives over $S$. 
This is made by taking
the pseudoabelianisation of $C \C(S)$ and then inserting Tate objects and 
twists by them:

\definition{2.9 Definition} An {\it object} of $C \M (S)$ is a triple
$$(X, P, r)$$
also denoted $(X,P)(r)$, where $X$ is a {\it smooth}, not necessarily
connected, variety, together with a {\it projective} morphism $f: X \to S$,
$P\in \En_{C \C S}(X,X)$ a projector, and $r \in {\Bbb Z}$ is
an integer.

{\it Morphisms} in $C \M (S)$ are defined as
$$\Ho_{C \M S} \bigl( (X,P,r), (Y,Q,s)\bigr)=
Q \circ \bigl( \oplus C_{\di Y_\alpha +r-s} (X\times_S Y_\alpha)\bigr) 
\circ P$$
where $Y = \coprod Y_\alpha$ is the decomposition of $Y$ into its connected
components $Y_\alpha$. The composition of morphisms is by means of the same
formula used for composing correspondences.
\enddefinition

\definition{2.10 Remarks, terminology, notation, etc}

(1) For $C = CH$, the category $CH \M(S)$ is called the category of {\it Chow}
motives over $S$. If $C=A$ is cycles modulo homological equivalence, 
$A \M (S)$ is the category of {\it homological} motives over $S$.
We will only need $A \M (S)$ very briefly in \S 4 and \S 5.
At this time, we are not in a position of constructing the analogue
of Grothendieck motives $\M$: this will be done in \S 4.

(2) Denoting $\Cal V(S)$ the category of smooth varieties $X$, projective
over $S$, there are natural contravariant {\it cohomological} $h_S: 
{\Cal V(S)} \to CH \M (S)$ and covariant {\it homological} $h^\vee_S
: {\Cal V(S)} \to CH \M(S)$ functors.
As an object, $h_SX =X$ regarded as a Chow motive, and for a morphism
$f: X \to Y$ covering the identity of $S$, $h_S(f)=
cl \Gamma_f^t$ is the cycle class of the transpose
$\Gamma_f^t \subset Y \times_S X$ of the graph $\Gamma_f$ of $f$.
Similarly, $h^\vee_S X = \oplus X_\alpha(\di X_\alpha)$, where $X=\coprod
X_\alpha$ is the decomposition in connected components, and
$h_S^\vee (f)=cl \Gamma_f$.

(3) If $f:X \to S$ is the morphism to $S$, we will sometimes use the 
following alternative notation for the 
objects $h_S X$ and $h_S^\vee X$:
$$\align
h_S X & = \CR f_\ast Q_X \\
h^\vee_S X& = \CR f_\ast D_X \endalign$$
In the coming subsection, we will construct a realisation functor
$CH \M (S) \to \De^b_{cc} (S)$. The notation is meant to suggest, for
instance, that $h_S X=\CR f_\ast Q_X$ realises to $Rf_\ast Q_X$, and
is therefore a {\it C}how theoretic ``$Rf_\ast$''. This notation will 
be particularly useful in the statement of the decomposition theorem in \S 6.
Similar remarks apply to the dualising sheaf $D_X$.

(4) As already said, realisations will be constructed in the coming subsection.
\enddefinition

\subhead {} Realisations \endsubhead

This subsection is devoted to the proof of the following basic result:

\proclaim{2.11 Theorem} There is a natural (faithful) realisation functor
$$A\M (S) \to \De^b_{cc}(S)$$
Therefore, composing with the forgetful functor $CH \M (S) \to 
A \M (S)$, there is also a realisation functor 
$$CH \M(S) \to \De^b_{cc} (S)$$
\endproclaim 

\definition{2.12 Remark} We have no special notation for the realisation
functors. In the instances where we will need to emphasise it, we will
simply denote it $\r$. For instance $\r : CH \M (S) \to \De^b_{cc} (S)$
is the realisation functor.
\enddefinition

\demo{Proof} The idea is to send the object $f: X \to S$ to the sheaf
$Rf_\ast Q_X$. The actual proof is made of the following ingredients:

(1) By lemma 2.14 below, the category $\De^b_{cc} (S)$ is pseudoabelian.

(2) Let $X$, $Y$ be smooth and $p: X \to S$, $q: Y \to S$ be projective
morphisms. By lemma 2.15 below, there is a natural isomorphism
$$\varphi : \Ho_S \bigl(Rp_\ast Q_X[2r], Rq_\ast Q_Y[2s]\bigr)
@>\cong >>H^{BM}_{2\di Y +2r -2s} X \times_S Y$$

(3) By lemma 2.17, the morphism $\varphi$ is compatible with composition,
i.e. $\varphi (v \circ u)= \varphi (v) \bullet \varphi (u)$. As for the
Chow theory, the last quantity is defined as
$p_{XZ \ast}\delta^! (v \times u)$: $\delta$ carries a natural 
orientation class in Borel-Moore homology, compatible with the 
Chow theoretic orientation class, see \cite{BFM}.

This allows to define realisations by sending
$$(X, P, r) \to \Img \bigl(\varphi^{-1} cl P : Rf_\ast Q_X \to 
Rf_\ast Q_X\bigr) [2r]$$ 
where $cl P \in H^{BM}_\bullet X \times_S X$ is the Borel-Moore class
of $P$. \qed \enddemo

The construction of the realisation functors relies on 3 lemmas which we
shall now state and prove. These are technical routines in $\De^b_{cc} S$, 
but basically easy
and we advise to skip them on first reading: after all, the statement
is rather plausible. Of these, the hardest is 2.17, stating that
composition in $CH \M (S)$ is compatible with composition in $\De^b_{cc} S$.
After some standard yoga, this is seen to be equivalent to the statement
that, on a smooth variety $Y$, cup product of cohomology classes is compatible,
via the Poincar\'e duality isomorphism, with intersection of
Borel-Moore homology classes. The key point is finally dealt with 
in 2.18, which will use the topological
bivariant theory of \cite{FM} 
in an essential way, and will be shown 
after a quick reminder of \cite{FM}.

\definition{2.13 Warning} Working in $\De^b_{cc} (S)$, we make no attempts
to get the signs right: in what follows all formulas are to be understood
up to a $\pm$ sign. This policy is well established in the literature
\cite{KS} \cite{SaM}.
\enddefinition 

\proclaim{2.14 Lemma} The category $\De^b_{cc}(S)$ is pseudoabelian.
\endproclaim

\demo{Proof} We will give a brief outline of a proof of this statement, for
which we could find no reference in the literature. In fact the proof
works for the full subcategory
$D^b$ of cohomologically bounded objects in any triangulated category
$D$ with $t$-structure. 

{\smc Step 1.} Let $p^2= p : M \to M$ be a projector, we wish to construct
the kernel and image $K$ and $I$ of $p$. 
The proof is by induction on the cohomological amplitude of $M$.
For a suitable $i$, $M^\prime = \tau_{\leq i} M$ and $M^{\prime \prime}
=\tau_{>i} M$ both have smaller cohomological amplitude then $M$, so
we may assume by induction that $p^\prime =\tau_{\leq i} p$, resp.
$p^{\prime \prime}=\tau_{>i} p$, have kernel and image $K^\prime$
and $I^\prime$, resp. $K^{\prime \prime}$ and $I^{\prime \prime}$.
Note that we have a morphism of exact triangles
$$\diagram
M^\prime \rto \dto^{p^\prime}&M \rto \dto^p&M^{\prime \prime}\rto^{[1]}
\dto^{p^{\prime \prime}}&\\
M^\prime \rto &M \rto &M^{\prime \prime}\rto^{[1]}&\enddiagram$$

{\smc Step 2.} With the identification $M^\prime = K^\prime \oplus
I^\prime$, $p^\prime: M^\prime \to M^\prime$ is the projection to 
the second factor, and similarly for $M^{\prime \prime}$.
Denoting $\delta: M^{\prime \prime} \to M^{\prime}[1]$ the
map of degree 1, we have a commutative diagram
$$\diagram 
K^{\prime \prime} \oplus I^{\prime \prime}\rto^\delta \dto
&K^{\prime}[1] \oplus I^{\prime}[1]\dto\\
K^{\prime \prime} \oplus I^{\prime \prime}\rto^\delta
&K^{\prime}[1] \oplus I^{\prime}[1]\enddiagram$$
where the vertical arrows are projection on the second factor.
We deduce that $\delta (K^{\prime \prime}) \subset K^\prime [1]$
(this has an obvious meaning in any additive category).
Similarly, arguing with $1-p$ instead of $p$, we also have
that $\delta (I^{\prime \prime}) \subset I^\prime [1]$.

{\smc Step 3.} Choose now a triangle
$$K^\prime \to K \to K^{\prime \prime}\to
K^\prime [1]$$
there is then a morphism $\varepsilon : K \to M$ so that the 
following is a morphism of triangles
$$\diagram
K^\prime\dto\rto &K\dto^\varepsilon \rto &K^{\prime \prime}\dto \rto^{[1]}& \\
M^\prime \rto &M \rto &M^{\prime \prime} \rto^{[1]}& \enddiagram$$
Replacing $\varepsilon$ with $\varepsilon -p\circ \varepsilon$, we
may assume that $p\circ \varepsilon =0$.

{\smc Step 4.} Let now $F : D \to$ abelian groups be any cohomological functor.
Then the following sequence is exact
$$0 \to F K \overset \varepsilon \to \to F M \overset p \to \to F M$$
Using $p\circ \varepsilon =0$ and step 2, the claim follows from
a never ending diagram chase along the paths and trails of the diagram
$$\diagram
0 \dto & 0 \dto &   & 0 \dto & 0 \dto \\
F^{-1}K^{\prime \prime}\rto\dto&FK^\prime\rto\dto &FK\rto\dto^\varepsilon &
FK^{\prime\prime}\rto\dto&F^1K^\prime\dto\\  
F^{-1}K^{\prime \prime}\oplus F^{-1}I^{\prime \prime}\rto\dto&
FK^\prime\oplus FI^\prime \rto\dto&
FM\rto\dto^p &
FK^{\prime\prime}\oplus FI^{\prime\prime}\rto \dto&
F^1K^\prime \oplus F^1I^\prime\dto\\  
F^{-1}K^{\prime \prime}\oplus F^{-1}I^{\prime \prime}\rto&
FK^\prime\oplus FI^\prime \rto&
FM\rto &
FK^{\prime\prime}\oplus FI^{\prime\prime}\rto&
F^1K^\prime \oplus F^1 I^\prime\enddiagram$$

{\smc Step 5.} Apply step 4 to $F=\Ho (U, -)$ where $U$ is an arbitrary object.
This shows that $K=\Ker (p)$. Then $I= \Ker (1-p)$.\qed \enddemo

\proclaim{2.15 Lemma} Let $p: X \to S$, $q: Y \to S$ be morphisms
of varieties, and consider the following fibre square
diagram

$$\diagram
  & X\times_S Y  \dlto_{q^\prime} \ddto^f \drto^{p^\prime} & \\
X \drto_p & &Y \dlto^q \\
          &S& \enddiagram$$
Then:

(1) For sheaves $F \in \De^b_{cc} X$, $G \in \De^b_{cc}Y$, there is a
natural isomorphism
$$Rf_\ast R{\Cal Hom}_{X\times_S Y}
(q^{\prime \ast} F, p^{\prime !} G)= R{\Cal Hom}_S (Rp_! F, Rq_\ast G)$$

(2) In particular, if $p$ is proper and $Y$ is smooth, there is a natural
isomorphism
$$\varphi:\operatorname{Hom}_{ S} 
\bigl(Rp_\ast {Q}_X[i], Rq_\ast {Q}_Y[j]\bigr) @>\cong >>
H^{BM}_{2\di Y +i-j}X\times_S Y$$
\endproclaim

\demo{Proof} (1) follows from Verdier duality and proper base change
$$\align 
Rf_\ast R{\Cal Hom}_{X\times_S Y} (q^{\prime \ast} F, p^{\prime !} G)
&=\text{(standard\; duality)\;} 
Rp_\ast R{\Cal Hom}_{X} (F, Rq^\prime_\ast p^{\prime !} G) \\ 
&=\text{(proper\; base\; change)\;}
Rp_\ast R{\Cal Hom}_{X} (F, p^!Rq_\ast G) \\
&=\text{(Verdier\; duality)\;}
R{\Cal Hom}_S (Rp_! F, Rq_\ast G)\endalign$$ 
If $p$ is proper and $Y$ is smooth, $Rp_!=
Rp_\ast$ and $D_Y=\Q_Y[2\di Y]$. Setting $F=\Q_X[i]$, $G=\Q_Y[j]$ in (1), 
and 
taking $H^0$, we obtain
$$\align
\Ho_{ S}\bigl(Rp_\ast {Q}_X[i], Rq_\ast {Q}_Y[j]\bigr)&=
\Ho_{ X\times_S Y}(\Q_{X \times_S Y}[i], p^{\prime !}\Q_Y[j]) \\
&=\Ho_{ X\times_S Y}(\Q_{X \times_S Y}[i], 
p^{\prime !}D_Y[-2\di Y+j]) \\
&=\Ho_{ X\times_S Y}(\Q_{X \times_S Y}, 
D_{X\times_S Y}[-2\di Y-i+j]) \\
&=H^{BM}_{2\di Y+i-j}X\times_S Y \endalign$$ 
that is, (2).
\qed \enddemo

\definition{2.16 Remark} In the proof of (1) we went 
from $X\times_S Y$ to $S$ passing through $X$. We could have gone there 
passing through $Y$, 
getting the same isomorphism.
\enddefinition

\proclaim{2.17 Lemma} Let $p_1:X \to S$, $p_2:Y \to S$, $p_3:Z \to S$ be 
morphisms of varieties with $p_1$, $p_2$ proper and $Y$, $Z$
smooth. Let $u: Rp_{1 \ast} \Q_X[i] \to Rp_{2 \ast} \Q_Y[j]$, $v:
Rp_{2 \ast} \Q_Y[j] \to Rp_{3 \ast} \Q_Z[k]$ be morphisms, and
$\varphi$ be the isomorphism of 2.15(2) above. Then
$$\varphi (v \circ u)=\varphi (v) \bullet \varphi (u)$$
\endproclaim 

The proof of 2.17 uses the following statement, with 
$Y$, $X\times_S Y$, $Y \times_S Z$ in place of $T$, $U$, $V$, and we 
postpone it until the end of the subsection.

\proclaim{2.18 Lemma} Let $T$ be a {\it smooth} variety, and
let $p: U \to T$, $q: V \to T$ be morphisms, with $p$ proper.
There are natural isomorphisms
$$\lambda^\prime :\Ho_{T}(Rp_\ast \Q_U ,\Q_T[i]) 
@>\cong >>H^{BM}_{2\di T-i}U$$
$$\mu^\prime : \Ho_{T}(\Q_T,Rq_\ast q^!\Q_T[j]) 
@>\cong >> H^{BM}_{2\di T-j}V$$
$$\nu^\prime : \Ho_{T}(Rp_\ast \Q_U,Rq_\ast q^!\Q_T[i+j]) 
@> \cong >> H^{BM}_{2\di T-i-j}U\times_T V$$
satisfying the identity
$$\nu^\prime (v \circ u)=\delta^!\bigl(
\mu^\prime (v)\times \lambda^\prime (u)\bigr)$$\endproclaim

The proof of 2.18 will keep us busy for some time. Let us explain the idea,
which is pretty basic. In the simplest case where both $p:U \to T$ and 
$q: V \to T$ are the identity map $T=T$, the lemma just says that cup
product in cohomology is compatible, via the Poincar\'e duality isomorphism,
with intersection product in homology. Indeed $\Ho_T (Q_T, Q_T[i])=
H^i T$, $\Ho_T (Q_T [i], Q_T[i+j])=H^{j} T$, and $v \circ u
\in \Ho_T(Q_T, Q_T[i+j])$ {\it is} the cup product $v \cup u$. Here
we let $\lambda^\prime = \mu^\prime = \nu^\prime =P: H^\bullet T
\to H^{BM}_{2\di T -\bullet} T$ be the Poincar\'e duality isomorphism.
The lemma then says $P (v \cup u) = \delta^! (Pv \times Pu)$, but this
is fine because $\delta^! (Pv \times Pu)=Pv \cdot Pu$ is the intersection
product. The actual proof of 2.18 uses the topological bivariant theory in
an essential way. We will now give a quick reminder of \cite{FM}, followed
by a generalisation of cup and intersection products, and finally the proof
of 2.18. In closing this subsection, we shall prove 2.17 using 2.18.

\subsubhead Topological bivariant theory \endsubsubhead

This is a very quick summary of the relevant bits of \cite{FM}. For more
information, the reader is invited to consult the original source.

(1) The topological bivariant theory associates to a morphism
$f: X \to Y$ of algebraic varieties $Q$-vector spaces $H^i (X \to Y)$.
Particular cases are $H^i (X =X)=H^i X$ is ordinary cohomology,
$H^i(X \to pt)=H^{BM}_{-i} X$ is Borel-Moore homology, and, for
the inclusion $X \hookrightarrow Y$ of a locally closed subvariety,
$H^i (X \hookrightarrow Y)=H^i_X Y$ is local cohomology.

(2) The natural operations are {\it products, proper push forward and 
pull back}.

If $\alpha \in H^i(X \to Y)$ and $\beta \in H^j(Y\to Z)$, there is a product
$\alpha \cdot \beta \in H^{i+j} (X \to Z)$.

If $f: X \to Y$ is proper and $Y \to Z$ is arbitrary, we
get a push forward homomorphism
$$f_\ast : H^i(X \to Z) \to H^i(Y \to Z)$$

If
$$\diagram 
X^\prime \dto \rto & X \dto \\
Y^\prime \rto^f    & Y \enddiagram$$
is a fibre square, we get a pull back
$$f^\ast : H^i(X \to Y) \to H^i(X^\prime \to Y^\prime)$$

(3) Proper push forward and pull back are functorial and satisfy
a number of natural compatibility axioms with products, like the 
projection formula, which we
will not bother writing down.

(4) A {\it strong orientation} for $f: X \to Y$ is a class $\omega \in H^j
(X \to Y)$ such that 
$$H^i(U \to X)\ni \alpha \to \alpha \cdot \omega
\in H^{i+j} (U \to Y)$$
is an isomorphism for all $U \to X$. With the aid of a strong
orientation we can define unexpected proper push forward
$f_! : H^\bullet X \to H^\bullet Y$ and pull back $f^! : H^{BM}_\bullet Y
\to H^{BM}_\bullet X$.
 
A class of maps closed under composition
possesses {\it canonical orientations} if all maps in the class are oriented
in a way that products of orientations are compatible with compositions
in the class. The key example of such a class is morphisms of smooth 
varieties (see below),
but also local complete intersection (lci) morphisms.

(5) One way to construct the topological bivariant theory is to do so on top
of the derived category $\De^b_{cc}$ as follows. For a morphism
$f: X \to Y$ one defines
$$H^i (X \to Y)=\Ho_Y(Rf_! Q_X, Q_Y[i])=\;\text{(Verdier duality)}\;
\Ho_X(Q_X, f^!Q_Y[i])$$

The {\it product} is essentially given by composition. Indeed let $f: X \to Y$,
$g:Y \to Z$ and $\alpha :Rf_!Q_X \to Q_Y[i]$, $\beta:Rg_!Q_Y \to Q_Z[j]$
be bivariant classes. The product $\alpha \cdot \beta$ is the composition
$$R(g\circ f)_! Q_X @> Rg_! \alpha >>
Rg_! Q_Y[i] \overset \beta \to \to Q_Z[i+j]$$
Note that we might as well have done the composition on $Y$ or on $X$,
with the same output.

If $f:X \to Y$ is proper, $g:Y \to Z$ is arbitrary, and
$\alpha :R(g\circ f)_! Q_X \to Q_Z[i]$ a bivariant class, the {\it proper
push forward} $f_\ast \alpha :Rg_! Q_Y \to Q_Z[i]$ uses the
canonical trace map $Q_Y \to Rf_\ast Q_X$ and is definded to be the composition
$$Rg_! Q_Y @> Rg_! tr >> Rg_!Rf_\ast Q_X=R(g\circ f)_!Q_X 
\overset \alpha \to \to Q_Z[i]$$

Finally, given a fibre square
$$\diagram 
X^\prime \dto_{g^\prime}\rto^{f^\prime}& X \dto^g \\
Y^\prime \rto_f    & Y \enddiagram$$
and a class $\alpha: Rg_!Q_X \to Q_Y[i]$, the {\it pull back} $f^\ast \alpha$
uses the base change isomorphism $Rg^\prime_! f^{\prime \ast}=
f^\ast Rg_!$:
$$Rg^\prime_! Q_{X^\prime}=
Rg^\prime_! f^{\prime \ast} Q_X =f^\ast Rg_!Q_X @>f^\ast \alpha >>
f^\ast Q_Y[i]=Q_{Y^\prime}[i]$$

If $X$ and $Y$ are smooth 
(and, for simplicity, equidimensional), $f: X \to Y$ possesses a {\it canonical
orientation}. Indeed, the duality homomorphism $Rf_! D_X \to D_Y$ is
here nothing but a morphism $Rf_! Q_X[2\di X] \to Q_Y[2\di Y]$, i.e.
a class in $H^{2\di Y -2 \di X} (X \to Y)$.

\subsubhead Cup and intersection products \endsubsubhead

\definition{2.19 Definition} Given $p:U \to T$ and $q: V\to T$ we define
a cup product
$$\cup : H^j(V \to T) \times H^i(U \to T) \to H^{i+j}(U\times_T V\to T)$$
by the formula
$$\alpha \cup \beta =q^\ast (\beta) \cdot \alpha$$
via the diagram
$$\diagram
U \times_T V \rto \dto & V \dto^q \rto & T \\
U \rto_p                 &T              & \enddiagram$$
(it is the same as $(-1)^{i+j}p^\ast \alpha \cdot \beta$).
\enddefinition

If $T$ is smooth, then it has a canonical orientation class
$\omega \in H^{-2\di T} (T \to pt)$. For any $W \to T$ we denote
$P_\omega$ the associated ``Poincar\'e duality'' isomorphism
$$H^i(W \to T) \ni \alpha @> P_\omega >> \alpha \cdot \omega
\in H^{i-2\di T} (W \to pt)$$
The diagonal morphism $\delta : T \to T \times T$ also has a canonical
orientation
$$or_\delta = P_{\omega \times \omega}^{-1} (\omega) \in H^{2\di T} (T
\to T \times T)$$
this allows 2 a priori different, but in the end equal, ways to define 
intersection products in Borel-Moore homology.

\definition{2.20 Definition-Proposition} Given $U \to T$ and $q: V\to T$ 
we define an intersection product
$$\bullet : H_i^{BM}U \times H_j^{BM}V \to H_{2\di T-i-j}^{BM} U\times_T V$$
in any of the 2 equivalent ways
$$a \bullet b= P_\omega (P_\omega^{-1} a \cup P_\omega^{-1} b) =
\delta^! (a \times b)$$
\enddefinition

\demo{Proof} To prove that 
$$P_\omega (P_\omega^{-1} a \cup P_\omega^{-1} b) =
\delta^! (a \times b)$$
the reader is invited to stare at the following diagram, where every square
is a fibre square
$$\diagram
      &   &    &U\dlto\xto'[1,0][2,0]&    &U\times_TV\dlto\ddto\llto \\
U\ddto&   &U\times T\llto\ddto&      &U\times V\llto\ddto&      \\
      &   &                   &T\dlto_\delta&        &V\dlto\xto'[0,-1][0,-2]\\
T\ddto&   &T\times T\llto\ddto&      &T\times V\llto\ddto&      \\
      &   &                   &      &                   &      \\
pt    &   &T\llto             &      &V\llto             &  \enddiagram$$
\qed \enddemo

\subsubhead Proof of 2.18 \endsubsubhead

{\smc Step 1.} We begin constructing natural isomorphisms
$$\lambda: \Ho_{ T}(Rp_\ast \Q_U ,\Q_T[i]) 
@> \cong >>H^i(U\to T)$$
$$\mu: \Ho_{ T}(\Q_T,Rq_\ast q^!\Q_T[j]) 
@> \cong >> H^j(V\to T)$$
$$\nu: \Ho_{ T}(Rp_\ast \Q_U,Rq_\ast q^!\Q_T[i+j]) 
@> \cong >> H^{i+j}(U\times_T V \to T)$$
satisfying the identity
$$\nu(v \circ u)=\mu (v)\cup \lambda (u)$$
The three isomorphisms are defined as follows. 
$\lambda$ is the identity
since ($p$ is proper) by definition
$$\Ho_{ T}(Rp_! \Q_U ,\Q_T[i]) = H^i(U\to T)$$
$\mu$ is defined with a single application of standard duality
$$\Ho_{ T}(\Q_T,Rq_\ast q^!\Q_T[j]) = 
\Ho_{V} (\Q_V, q^!\Q_T[j])= H^j(V \to T)$$
For the rest of the proof, we fix the notation in the following diagram
$$\diagram
  & U\times_T V  \dlto_{q^\prime}  \drto^{p^\prime} & \\
U \drto_p & &V \dlto^q \\
          &T& \enddiagram$$
To define $\nu$, let first $$\chi : \Ho_{V}(q^\ast Rp_\ast
Q_U, q^! Q_T[k])@>\cong >>  H^k(U\times_T V \to T)$$
be the isomorphism obtained composing the following natural identifications
$$\align
 &\Ho_{ V}
(q^\ast Rp_\ast \Q_U, q^!\Q_T[k]) \\
= (\text{base\; change},\; p\;\text{proper})\;
&\Ho_{ V}(Rp^\prime_\ast q^{\prime \ast} \Q_U, q^!\Q_T[k]) \\
=&\Ho_{ V}(Rp^\prime_\ast \Q_{U\times_T V}, q^!\Q_T[k]) \\
=(p\; \text{proper})\; &\Ho_{ U\times_T V}(
\Q_{U\times_T V}, p^{\prime !} q^! \Q_T[k]) \\
=&H^{k}(U\times_T V \to T)\endalign$$
As a small digression, let $u: Rp_\ast Q_U \to Q_T[i]$, $v: 
Q_T \to Rq_\ast q^! Q_T[j]$, and let $v^\prime : q^\ast Q_T=Q_V \to
q^! Q_T [j]$ correspond to $v$ under Verdier duality. We like to observe, 
at this point, that
$$\mu (v) \cup \lambda (u)
= q^\ast\bigl(\lambda(u)\bigr)\cdot \mu (v)=\chi\bigl(v^\prime 
\circ q^\ast (u) \bigr)$$ 
This ends the digression. Now, to come back to the definition 
of $\nu$, we just compose $\chi$ 
with a standard duality isomorphism
$$\align
\Ho_{T}(Rp_\ast Q_U, Rq_\ast q^! Q_T[k])=
&\Ho_{V}(q^\ast Rp_\ast Q_U, q^! Q_T[k])\\
@> \chi >>&H^{k}(U\times_T V \to T)\endalign$$
Look now at the commutative diagram
$$\diagram
Rp_\ast\Q_U \dto_u\rto^{tr}&Rq_\ast q^\ast 
Rp_\ast\Q_U \dto^{Rq_\ast q^\ast \lambda (u)}\\
\Q_T[i]\rto^{tr}\dto_v&Rq_\ast q^\ast \Q_T[i]
\dto^{Rq_\ast v^\prime}\\
Rq_\ast q^! \Q_T[i+j]\rdouble &Rq_\ast q^! \Q_T[i+j] \enddiagram$$
The diagram shows that $v^\prime \circ q^\ast u$ corresponds to
$v \circ u$ under the standard duality isomorphism
$$\Ho_{ T}(Rp_\ast Q_U, Rq_\ast q^! Q_T[i+j])=
\Ho_{ V}(q^\ast Rp_\ast Q_U, q^! Q_T[i+j])$$
which was used in the definition of $\nu$. Therefore
$$\nu(v \circ u) = \chi (v^\prime \circ q^\ast u)$$
On the other hand, as we have seen in the digression
$$\mu (v) \cup \lambda (u)= \chi (v^\prime \circ q^\ast u)$$
and combining the last 2 displayed formulas concludes step 1.

{\smc Step 2.} We now define
$$\lambda^\prime = P_\omega \circ \lambda: 
\Ho_{T}(Rp_\ast \Q_U ,\Q_T[i]) 
@> \cong >> H^{BM}_{2\di T-i}U$$
$$\mu^\prime = P_\omega \circ \mu: \Ho_{T}(\Q_T,Rq_\ast q^!\Q_T[j]) 
@> \cong >> H^{BM}_{2\di T-j}V$$
$$\nu^\prime = P_\omega \circ \nu: \Ho_{T}(Rp_\ast \Q_U,Rq_\ast q^!\Q_T[i+j]) 
@> \cong >> H^{BM}_{2\di T-i-j}U\times_T V$$
The statement now follows from step 1, 2.20, and a simple calculation
$$\align
\nu^\prime (v \circ u)&=P_\omega \bigl(\nu (v\circ u)\bigr)=P_\omega
\bigl(\mu (v)\cup \lambda (u)\bigr)=\\
&=P_\omega\bigl(P_\omega^{-1}\mu^\prime(v)\cup P_\omega^{-1}\lambda^\prime(u)
\bigr)=\delta^!\bigl( 
\mu^\prime (v)\times\lambda^\prime (u) \bigr)\endalign$$
This finishes the proof of 2.18.

\subsubhead Proof of 2.17 \endsubsubhead

We summarise the notation for the various spaces and maps in 
the following commutative diagram
$$\diagram
           &X\times_S Z                                    &             \\
           &X\times_S Y \times_S Z\dlto\drto \uto_{p_{XZ}} &             \\
X\times_S Y\dto^{p_2^\prime}\drto^{p_1^\prime}& & Y\times_S Z
\dto^{p_2^{\prime\prime}}\dlto_{p_3^\prime} \\
X\drto_{p_1}&    Y\dto^{p_2}        &Z\dlto^{p_3}\\
            &    S                  &             \enddiagram$$

The proof of the lemma results from contemplating the following commutative 
diagram, which is commented upon below
$$\diagram
Rp_{1\ast} \Q_X[i]\dto^u\rto^{}&Rp_{2\ast}p_2^\ast Rp_{1\ast}\Q_X[i]
\dto^{Rp_{2\ast}u^\prime}\rdouble^{}&Rp_{2\ast}Rp^\prime_{1\ast}
\Q_{X\times_S Y}[i]\dto^{Rp_{2\ast}u^\prime} \\
Rp_{2\ast} \Q_Y[j]\dto^v\rdouble &Rp_{2\ast}\Q_Y[j]
\dto^{Rp_{2\ast}v^\prime}\rdouble       &Rp_{2\ast}\Q_Y[j]
\dto^{Rp_{2\ast}v^\prime} \\
Rp_{3\ast} \Q_Z[k]     &Rp_{2\ast}p_2^!Rp_{3\ast}\Q_Z[k]\lto^{}
\rdouble^{}                         &Rp_{2\ast} Rp^\prime_{3\ast}
p^{\prime !}_3 \Q_Y[2\di Y-2\di Z +k] \enddiagram$$
Here $u^\prime :p_2^\ast Rp_{1\ast}\Q_X[i]\to \Q_Y[j]$ corresponds to $u$ 
via the standard duality 
$$\Ho_Y(p_2^\ast Rp_{1\ast}\Q_X[i], \Q_Y[j])=\Ho_S (Rp_{1\ast} \Q_X[i],
Rp_{2\ast} \Q_Y[j])$$
and $tr: Rp_{1\ast} \Q_X[i] \to Rp_{2\ast}p_2^\ast Rp_{1\ast}\Q_X[i]$
is the trace map giving rise to 
(or arising from, depending on the reader's preference) the duality. 
Similar comments apply to 
$v^\prime$. 

The two equal signs in the right portion of the diagram
are obtained from the proper base change isomorphism, 
for instance the lower one (which is the hardest) is derived as follows
$$\align
p_2^!Rp_{3\ast}\Q_Z[k]= \text{(base\; change)\;} &Rp^\prime_{3\ast}
p^{\prime \prime !}_2 \Q_Z[k] \\
= (Z\; \text{is \;smooth})\; &Rp^\prime_{3\ast}
p^{\prime \prime !}_2 D_Z[-2\di Z+k]  \\
= &Rp^\prime_{3\ast}D_{Y\times_S Z}
[-2\di Z+k]  \\
=(Y \; \text{is \;smooth})\; &Rp^\prime_{3\ast}
p^{\prime !}_3 \Q_Y[2\di Y-2\di Z +k] \endalign$$
Contemplating the diagram in the light of how $\varphi$ is defined (lemma 2.15 
and remark 2.16), 
and using 2.18, 
with $X\times_S Y$, $Y\times_S Z$
and $Y$ in place of $U$, $V$ and $T$ respectively, we evince the 
following
$$\varphi (v\circ u)=p_{XZ \ast} \nu^\prime (v^\prime \circ u^\prime)$$
$$\varphi (u)=\lambda^\prime (u^\prime)$$
$$\varphi (v)=\mu^\prime (v^\prime)$$
The result then follows immediately from 2.18. This finishes the proof of 2.17.

\head {3} Standard conjectures and canonical filtrations \endhead

In this section, which is intended mainly for reference, we begin recalling
Gro\-then\-dieck's standard conjectures, which were introduced, 
among other things,
to determine the behaviour of the category $\M$ of Grothendieck motives.
We will mainly need them in \S 4, when we will define the relative analogue
$\M (S)$ of $\M$ and show, assuming the conjectures, that it is an abelian
semisimple category and the decomposition theorem holds in $\M(S)$. Then
we recall Murre's conjecture, which we only need later on in \S 5 and
\S 6, implying the existence of a natural
filtration $F^\bullet$ on the Chow groups of smooth and projective varieties, 
and
explain how this conjecture can be used to fill in part of
the gap between $\M$ and $CH \M$, making it possible to define
a noncanonical decomposition of a Chow motive into its cohomology
groups.
Finally, we recall S. Saito's unconditional definition of a filtration,
having all the expected categorical properties, except that it is not
known to be separated.
If it is separated, it coincides with Murre's, and Murre's
conjecture holds. 

It is Saito's filtration that will be extended, in \S 5,
to the Chow groups of arbitrary quasiprojective varieties.
This will be used in the final \S 6 to prove the decomposition theorem
in $CH \M (S)$, and in the forthcoming part II, when we will
propose an unconditional definition of the intersection Chow groups
$ICH^r X$.

\subhead {} Standard conjectures \endsubhead

Before stating the conjectures, we introduce some notation and recall some 
well known facts on the cohomology of algebraic varieties. 

Let $X$ be a smooth projective variety of complex dimension $d$, with a fixed
ample divisor class $L \in H^2 X$.
The Lefschetz theorem asserts that, for $i \leq d$, the $d-i$-th iterated
cup product with $L$ is an isomorphism of $H^iX$ to $H^{2d-i}X$
$$L^{d-i}: H^iX @> \cong >> H^{2d-i}X$$
For $i \leq d$ we then define the primitive cohomology of $X$ to be
$$P^iX= \Ker L^{d-i+1} \subset H^i X$$
We have the hard Lefschetz decomposition of the 
cohomology of $X$
$$H^iX=\oplus_{j \geq 0} L^jP^{i-2j}X$$
if $i\leq d$ and
$$H^iX=\oplus_{j\geq i-d}L^jP^{i-2j}X$$
if $i>d$.

\definition{3.1 Definition} The Lefschetz operator $\Lambda: H^i X \to 
H^{i-2}X$ relative to the ample class $L$ 
is defined as follows. Let $\alpha \in H^i X$ and write, using the
Lefschetz decomposition 
$$\alpha = \sum_j L^j \alpha^{i-2j}$$
with $\alpha^{i-2j}\in P^{i-2j}X$. Then by definition
$$\Lambda \alpha =\sum_j L^{j-1} \alpha^{i-2j}$$
(i.e., $\Lambda$ removes one $L$). 
\enddefinition

Grothendieck \cite{Gr} proposed the following 2 {\it standard conjectures}:

\proclaim{3.2 Standard conjecture of Lefschetz type} 
The $\Lambda$ operator is algebraic.
\endproclaim

\proclaim{3.3 Standard conjecture of Hodge type} The rational quadratic 
form
$$(\alpha, \beta) \to (-1)^i tr(\alpha \cup \beta \cup L^{d-2i})$$
is positive definite on $P^{2i} \cap H^{2i}X_{alg}$. 
\endproclaim

Progress is occasionally made on these conjectures \cite{Ja1} \cite{Sm}.

In the proof of the decomposition theorem in $\M (S)$, \S 4, we will need 
the following simple consequence of the conjecture of Lefschetz type:

\proclaim{3.4 Proposition} Assume the standard conjecture of Lefschetz type. 
Let $S$ be
a smooth quasiprojective variety, and
$f:X \to S$ be a smooth projective morphism, with relatively
ample divisor class $L\in H^2 X$. There exists a cycle
$$Z \in CH_{\di X+1} (X \times_S X)$$
such that, for every $s \in S$ and fibre $X_s$, 
$Z|_{X_s \times X_s}$ induces the $\Lambda_s$ operator (relative
to the class $L_s=L|{X_s}$) of that fibre.
\endproclaim
 
\demo{Proof} The proof uses a standard ``spreading out'' argument
followed by specialisation.
By conjecture 1, there is a cycle 
$Z_{\eta}$ on $X_{\eta}
\times X_{\eta}$ inducing  
$\Lambda_{\eta}$. 
Let
$U \subset S$ be a neighbourhood of $\eta$ and $Z_U$ a cycle on $X_U\times_U
X_U$ such that $Z_U|\eta=Z_\eta$,
and let $Z$ on $X\times_S X$ be its Zariski closure. I claim that for all
scheme theoretic points $s \in S$, 
$Z|_{X_s\times X_s}$
induces $\Lambda_s$. By considering a chain of points
$$s \in \overline {s_1} \in \overline {s_2} \in \cdots \in
\overline {\eta}$$
with 
$$\operatorname{cod}_{s_i} \overline {s_{i+1}}=1$$
we are reduced to the case of
the spectrum $T$ of a discrete valuation 
ring with central point $0$ and generic point $s$, and a morphism
$T \to S$.
Assuming that $Z|_{X_s\times X_s}$ induces $\Lambda_s$, we 
need to prove that
$Z|_{X_0\times X_0}$ induces $\Lambda_0$. 
In this situation, letting $k(s)$ be the function field of $T$,
there are well defined specialisation maps
$$\diagram
CH^i X_s\times X_s \rto^{csp} \dto_{cl} & CH^i X_0\times X_0 \dto^{cl}\\
H^{2i} X_s \times X_s \rdouble_{hsp}   & H^{2i}X_0 \times X_0 \enddiagram$$
Let us recall the construction of the Chow theoretic specialisation 
homomorphism. We have a diagram
$$\diagram
CH^{i-1} X_0\times X_0 \rto^{i_\ast}&CH^iX\times_TX\rto^{j^\ast}
\dto^{i^!}&CH^iX_s\times X_s\rto&0\\
          &CH^i X_0 \times X_0&  &  \enddiagram$$
where the row is exact. Because $X_0 \sim 0$, we can define
$csp (\alpha) = i^! \alpha^\prime$ where $\alpha^\prime \in 
CH^iX\times_TX$ is anything such that $\alpha^\prime|_{X_s \times X_s} 
=\alpha$, 
and the result does not depend on $\alpha^\prime$. By construction of
$csp$ then 
$$Z|_{X_0\times X_0}= csp (Z|_{X_s\times X_s})$$
On the other hand clearly $\Lambda_0 =hsp (\Lambda_s)$,
indeed by what we just said about specialisation $L_0=csp (L_s)$, so ``removing
one $L$'' specialises to ``removing one $L$''. Therefore
$$cl (Z|_{X_0\times X_0}) =\Lambda_0$$ \qed \enddemo

The most important consequence of the standard conjectures is the following 
\cite{Kl1, Kl2}:

\proclaim{3.5 Theorem} Assuming the standard conjectures, then

(1) $A \M=\M$, in other words homological and numerical equivalence
of algebraic cycles are the same.

(2) The category
$\Cal M $ of Grothendieck motives is abelian and semisimple.
\endproclaim

\demo{Proof} See \cite{Kl1}. \qed\enddemo

\subhead {} Murre's conjecture
\endsubhead 

The standard conjectures are perfectly adequate in determining the behaviour
of Grothendieck motives. There is a large gap between Grothendieck motives
and Chow motives, which one begins to appreciate when trying to decompose
a Chow motive $hX$ into its pieces $h^iX[-i]$ in an unambiguous way.
To address this issue, Murre \cite{Mu} proposed the following:

\proclaim{3.6 Murre's conjecture} Let $X$ be a smooth variety of complex 
dimension $d$,
and $\pi^i \in H^iX \otimes H^{2d-i}X \subset H^{2d}X\times X$ be the 
K\"unneth components of the diagonal. Then:

(A) The $\pi^i$ lift to an orthogonal set of projectors $\Pi^i \in 
CH^{\di X} (X \times X)$ such that $\Delta = \sum \Pi^i$.

(B) The correspondences $\Pi^{2r+1}$,...,$\Pi^{2\di X}$ act as zero on
$CH^r X$.  

(C) For each $\nu$, $F^\nu CH^rX= \Ker\Pi^{2r}\cap \dots \cap 
\Ker \Pi^{2r-\nu+1}$
is independent of the choice of the $\Pi^i$.

(D) $F^1CH^r X=CH^r X_{hom}$ is the group of cycles homologically equivalent 
to zero.
\endproclaim

Jannsen proved \cite{Ja2, pg. 294--296 and 259}:

\proclaim{3.7 Theorem} Assume 3.6. The filtration $F^\bullet$ satisfies the
following properties:

(a) $F^0 CH^rX=CH^rX$, $F^1CH^r X=CH^r X_{hom}$.

(b) $F^\nu CH^r X \cdot F^\mu CH^s X \subset F^{\nu+\mu}CH^{r+s}X$.

(c) If $f:X \to Y$ is a morphism of smooth projective varieties, $f_\ast$ 
and $f^\ast$ respect the filtration (no shifts involved). 

(d) Let $\Gamma \in CH_{\di X}X \times X$ be a correspondence. Assume that
$\Gamma$ acts trivially on $H^{2r-\nu}X$. Then
$$\Gamma_\ast: F^\nu CH^r X \to F^{\nu + 1} CH^r X$$

(e) Assuming moreover the standard conjecture of 
Lefschetz type, $F^{r+1} CH^r X =(0)$. \qed
\endproclaim

We will show momentarily that, assuming the standard conjecture
of Lefschetz type, Murre's $A+B+C+D$ is equivalent to $A+B^\prime +D$
below:

\proclaim{3.8 Conjecture $B^\prime$ (vanishing)} Let $X$ be a smooth 
projective variety
and $P \in CH^{\di X} X\times X$ a projector. Assume $P_\ast H^iX=0$
for $i\leq 2r$. Then $P_\ast CH^r X=(0)$.
\endproclaim

\proclaim{3.9 Proposition} Assuming the standard conjecture
of Lefschetz type, Murre's $A+B+C+D$ implies $B^\prime$.
\endproclaim

\demo{Proof} Assume $A+B+C+D$. By 3.7(d), 
$P_\ast F^\nu CH^r X \subset F^{\nu+1} CH^r X$ for
all $\nu \leq 2r$. But $P$ is a projector, so
$P_\ast CH^rX =P^2_\ast CH^r X \subset P_\ast F^1 CH^r X= 
P^2_\ast F^1CH^r X \subset P_\ast F^2 CH^r X... \subset F^{r+1}CH^rX=0$ 
by 3.7(e). \qed \enddemo

To prove that $A+B^\prime + D$ implies $A+B+C+D$ we now make a small 
digression to discuss the 
decomposition of Chow motives. The decomposition theorem in $CH\M (S)$
in \S 6 will follow the same basic strategy.

\definition{3.10 Definition} A Chow motive $M$ has cohomological degree
$\leq m$, resp. $\geq m$ if the cohomology groups
$$H^i M=0$$
vanish for $i>m$, resp. $i<m$. $M$ has degree exacly $m$ if it has degree
$\geq m$ and $\leq m$.
\enddefinition

\proclaim{3.11 Proposition} Conjecture $B^\prime$ implies

(1) If $M$ has cohomological degree $\leq m$ and $N$ has cohomological
degree $>m$, then
$$\Hom_{CH \M}(M,N)=0$$

(2) If $M$ and $N$ have cohomological degree exactly $m$, 
then the natural homomorphism
$$\Hom_{CH \M} (M,N) \to \Hom_\M (M,N)$$
is an isomorphism.
\endproclaim

\demo{Proof} Assume for simplicity that $M=(X,P)$ and $N=(Y,Q)$. Let
$Z=X \times Y$ and consider the {\it projector}
$$CH^i Z \ni \alpha \to \Psi \alpha =Q\circ \alpha \circ P\in CH^i Z$$
To prove both statements, it is enough to show that $\Psi =0$ if
$cl \Psi =0$, but this is $B^\prime$.\qed \enddemo

We draw 2 consequences

\proclaim{3.12 Corollary (decomposition of Chow motives)} 
Assume $A+B^\prime$, then for all smooth projective varieties $X$

(1) There is a noncanonical direct sum decomposition
$$h X = \sum (X, \Pi^i)$$

(2) The monomorphisms
$$\tau_{\leq i} hX= \sum_{m\leq i} (X, \Pi^m)\to hX$$
(where the 1st equality is a definition of $\tau_{\leq i} hX$)
are specified up to canonical isomorphism. In particular so
are the ``subquotients''
$$h^i X[-i]=(X, \Pi^i)$$
(this is a definition of $h^i X[-i]$) specified up to canonical 
isomorphism.\qed\endproclaim

\proclaim{3.13 Corollary} $A+B^\prime+D$ implies $A+B+C+D$.
\endproclaim

\demo{Proof} We can see that
$$F^\nu CH^i X=\Ho_{CH \M}\bigl(pt, (\tau_{\leq 2i-\nu} hX)(i)\bigr)$$
is independent on the $\Pi^i$s, by 3.12. \qed \enddemo

\subhead Saito's filtration \endsubhead

S. Saito \cite{SaS} gave an unconditional definition of a filtration
on the Chow groups of smooth projective algebraic varieties over $k$ and
proved that, assuming the standard conjectures, it coincides with Murre's
filtration. We now recall Saito's definition and his results:

\definition{3.14 Definition} \cite{SaS} For a smooth projective variety $X$ 
we define a filtration
$$CH^rX =F^0CH^rX \supset F^1CH^rX\supset \cdots \supset F^\nu CH^rX \supset
\cdots$$
in the following inductive way:

(1) $F^0CH^rX=CH^rX$.

(2) Assume $F^\nu CH^r X$ defined for all $X$ and all $r$. Then we set:
$$F^{\nu +1}CH^rX=\sum_{Y, q, \Gamma} \Gamma_\ast F^\nu CH^{r-q} Y$$
where $Y$, $q$ and $\Gamma$ range over the following data:

(2.1) $Y$ is smooth and projective,

(2.2) $q$ is an integer (the operation yields nothing unless
$r -\di Y \leq q \leq r$, since otherwise $CH^{r-q} Y =0$),

(2.3) $\Gamma \in CH^{\di Y+q} (Y \times X)=\Ho_{\CH}\bigl(Y, X(q) \bigr)$ 
is a correspondence such that
$$\Gamma_\ast H^{2r-2q-\nu}Y \subset N^{r-\nu+1} H^{2r-\nu}X$$
where $N^\bullet$ is Grothendieck's coniveau filtration (see 3.15 below).
\enddefinition

\definition{3.15 Reminder} Recall that the coniveau filtration on the 
cohomology of a smooth projective algebraic variety $X$ is defined as 
$$N^pH^iX=\sum_{Y,f} f_\ast H^{i-2q}Y$$
where the sum ranges over all smooth projective $Y$ with $q=\di X -\di Y
\geq p$ and morphisms $f: Y \to X$.
\enddefinition

\proclaim{3.16 Theorem} \cite{SaS} 

(1) The filtration defined in 3.14 satisfies the 
properties b, c, d in 3.7 and
$F^1CH^rX=CH^rX_{hom}$.

(2) If the filtration is separated, i.e. $F^\nu CH^i X = 0$ for $\nu$ large,
then Murre's conjecture is true and the filtrations are the same.

(3) Assuming the standard and Murre's conjectures, the filtrations are the same
(in particular they are separated). \qed
\endproclaim

\definition{3.17 Remark} It would have been possible to define $F^\bullet$ just
as in 3.14, but replacing 3.14(2.3) with the easier condition
$$0=\Gamma_\ast H^{2r-2q-\nu}Y \to H^{2r-\nu}X$$
3.16 would still be true (in fact, slightly easier to prove) for this 
filtration.
\enddefinition 

The following property is not stated explicitly in \cite{SaS} and will
be used in \S 5.

\proclaim{3.18 Proposition} For $X$, $Y$ smooth projective and $f:X \to Y$ 
a morphism, 
$$f_\ast : CH_sX \to CH_s Y$$
is strictly compatible with the $F^\bullet$ filtration as defined in 3.14.
\endproclaim

\demo{Proof} Choose a diagram
$$\diagram
Z\rto^i \drto_\pi & X \dto^f\\
                  & Y \enddiagram$$
where $Z$ is smooth projective, $i: Z \hookrightarrow X$ is a closed 
embedding, and $\pi
: Z \to Y$ a generically finite morphism of degree $d$. Let $\alpha
\in F^\nu CH^r X$. Then 
$$\alpha = {1\over d} \pi_\ast \pi^\ast \alpha= {1\over d}f_\ast
i_\ast \pi^\ast \alpha$$
It is clear that
$$i_\ast \pi^\ast \alpha \in F^\nu CH^r X$$  
\qed \enddemo

\head 4. Grothendieck motives over a base, semisimplicity and 
decomposition \endhead

This section is divided into 3 subsections. In the first, for the convenience 
of the reader and to fix the notation, we recall the notion of perverse 
sheaves and the statement of the topological decomposition theorem. The 
standard references are \cite{Bo}, 
\cite{BBD}. For the expert, we say right away that we found it convenient 
to use Deligne's convention for perverse sheaves, because better suited 
for taking direct images under a closed embedding, and Borel's convention
for intersection complexes. With our conventions, therefore, intersection 
complexes are not perverse (but a suitable shift is). In the second subsection
we define a category $\M (S)$, which we call the category of 
Grothendieck motives over a variety
$S$: this is the correct analogue of the category of Grothendieck motives
over the point and is built precisely 
in order to have a faithful realisation
in the graded category of perverse sheaves. In the third subsection, assuming 
the standard conjectures, 
we prove that $\M (S)$ is abelian and semisimple 
and, as a consequence, we derive a decomposition theorem in $\M (S)$ which
realizes to the topological decomposition theorem.

\subhead {} Perverse sheaves and the topological
decomposition theorem \endsubhead

In this subsection, we recall the theory of perverse sheaves and the
topological decomposition theorem. The standard reference for this 
material is \cite{BBD}. For ease of notation and
terminology, we will assume that $k = {\Bbb C}$ and refer the reader
to the original source for the language suitable to the \'etale situation.

\definition{4.1 Definition} Let $\De$ be a triangulated category. A
$t$-{\it structure} on $\De$ is a pair $(\De^{\leq 0}, \De^{\geq 0})$ of full
subcategories of $\De$, satisfying the following axioms:

(1) $\De^{\leq 0}[1] \subset \De^{\leq 0}$ and
$\De^{\geq 0} \subset \De^{\geq 0}[1]$.

(2) $\Ho (\De^{\leq 0}, \De^{>0})=0$.

(3) For every object $K$ of $\De$, there is a (necessarily unique up to 
canonical isomorphism) 
triangle
$$K^\prime \to K \to K^{\prime \prime} \overset 1 \to \to$$
with $K^\prime \in \De^{\leq 0}$, $K^{\prime \prime} \in \De^{>0}$.
\enddefinition

The assignment $K$ to $K^\prime = \tau_{\leq 0} K$ is functorial and the
corresponding functor is called the {\it truncation} functor relative to
the $t$-structure. $\tau_{\leq m}K$ is defined to be 
$\bigl(\tau_{\leq 0}(K[m])\bigr)[-m]$, similarly $\tau_{\geq m}$, and
${\Cal H}^m(-)= \bigl(\tau_{\leq m} \tau_{\geq m}(-)\bigr)[m]$
is the $m$-th {\it cohomology} functor relative to the $t$-structure.

The main theorem \cite{BBD} about $t$-structures 
asserts that the {\it heart}~ $
\De^{\leq 0} \cap \De^{\geq 0}$ is an abelian category. 

We now come to the most important example of $t$-structure, the {\it
perverse} $t$-structure on $\De^b_{cc}(S)$, but first:

\definition{4.2 Definition} Let $S$ be a quasi projective variety over a field
$k$. A {\it good} stratification of $S$ is a stratification
$$S=\coprod T_k$$
where $T_k$ is a Zariski locally closed subset of complex dimension $k$,
satisfying the following axioms:

(1) each stratum $T_k$ is smooth,

(2) the stratification is topologically normally locally trivial.
\enddefinition

From now on, we will assume that all varieties $S$ 
are equipped with a good stratification.

\definition{4.3 Notation} If $S = \coprod_k T_k$ is a good
stratification, we denote $i_{T_k}: T_k\to S$ 
the inclusion and $S_k =\coprod_{h\leq k} T_h$ the Zariski closure 
of $S_k$.
\enddefinition

\definition{4.4 Definition} 

(1) Let ${\Cal T}=\{T_k\}$ be a good stratification 
$$\De_{\Cal T}(S)=\{K\in \De^b (S) \mid K|{T_k}\; 
\text{is cohomologically locally constant}\;\forall k\}$$ 

(2) The bounded derived category of cohomologically constructible
sheaves is defined as
$$\De^b_{cc}(S)= \cup \De_{\Cal T}(S)$$
the union being taken over all good stratifications of $S$. 
\enddefinition

From now on, when dealing with a sheaf $K \in \De^b_{cc}(S)$, in
connection with a preexisting good stratification 
${\Cal T}=\{T_k\}$, we will assume
that $K$ is cohomologically locally constant along all strata
$T_k$ of ${\Cal T}$.

\definition{4.5 Definition} The {\it perverse} $t$-structure on
$\De^b_{cc}(S)$ is defined as follows
$$\align
{}^p\De^{\leq 0}&=\cup{}^p\De_{\Cal T}{}^{\leq 0}\\
{}^p\De^{\geq 0}&=\cup{}^p\De_{\Cal T}{}^{\geq 0}\endalign$$
the union being taken over all good stratifications ${\Cal T}$, where
$$\align
{}^p\De_{\Cal T}{}^{\leq 0}&=\{K\in {\Cal D}_{\Cal T} 
\mid {\Cal H}^i i_{T_{k}}^\ast K =0, \; i>-k\}\\
{}^p\De_{\Cal T}{}^{\geq 0}&=\{K\in {\Cal D}_{\Cal T} 
\mid {\Cal H}^i i_{T_{k}}^! K =0, \; i<-k\}
\endalign$$
\enddefinition
It is well known that the above data define a $t$-structure, whose heart
$= \Pe (S)$ is the category of {\it perverse sheaves} on $S$.
We will denote the truncation, resp. cohomology functors of the perverse
$t$-structure with the symbol ${}^p \tau_{\leq 0}$, resp. ${}^p {\Cal
H}^m$.

The most important construction in the theory of perverse sheaves is that
of the intersection complexes:
\smallskip

\noindent{\bf 4.6 Intersection complexes.} Let ${\Cal T}$ be a good
stratification and $V$ a rational local system on the
largest stratum $T_d$.
The {\it intersection complex} $\IC V$ is characterised by the properties
$$\align
&\IC V|{T_d}=V  \tag 0\\
&{\Cal H}^i i_{T_{k}}^\ast \IC V =0, \; i\geq d-k \;(k <d) \tag --\\
&{\Cal H}^i i_{T_{k}}^!    \IC V =0, \; i\leq d-k \;(k <d) \tag +
\endalign$$
It follows immediately from the characterisation just given that
the shift $\IC V [d]$ is a perverse sheaf on $S$. 

It is important to
understand that $\IC V[d]$ is not characterised by being a perverse sheaf 
and restricting to $V[d]$ on the largest stratum.
In fact, there are lots and lots of perverse
sheaves which restrict to $V[d]$ on $T_d$, and $\IC V[d]$ is built to be
as much in the middle of $\Pe (S)$ as possible. 

The {\it intersection cohomology} of $S$ is defined as
$$IH^m S=H^m \IC {Q}_{T_d}$$

Similarly, if $V_k$ is a local system on the stratum $T_k$ of dimension
$k$, the intersection complex $\IC V_k$ is a 
complex supported on the Zariski closure $S_k$ characterised by the properties
$$\align
&\IC V_k|{T_k}=V_k  \tag 0\\
&{\Cal H}^i i_{T_{h}}^\ast \IC V =0, \; i\geq k-h \;( h <k) \tag --\\
&{\Cal H}^i i_{T_{h}}^!    \IC V =0, \; i\leq k-h \;( h <k)  \tag +
\endalign$$
It follows immediately from the characterisation just given that
$\IC V_k [k]$ is a perverse sheaf on $S$.

\proclaim{4.7 Topological decomposition theorem} Let $X$ be a smooth
variety and $f: X \to S$ be a projective morphism. 

(1) There is a noncanonical direct sum decomposition
$$Rf_\ast {Q}_X \cong \sum {}^pR^mf_\ast {Q}_X [-m]$$
in $\De^b_{cc}(S)$, where ${}^pR^mf_\ast {Q}_X$ denotes
the $m$th perverse cohomology of $Rf_\ast {Q}_X$. 
The decomposition itself is not unique, but the subobjects
$${}^p \tau_{\leq m} Rf_\ast {Q}_X=\sum_{i\leq m}
{}^pR^if_\ast {Q}_X [-i]$$
are uniquely specified.

(2) Let ${\Cal T}=\{ T_k\}$ be a good stratification with the property that 
${}^p R^m f_\ast {Q}_X \in \De_{\Cal T}(S)$. There are local systems
$V^{m}_k$ on $T_k$ and a canonical isomorphism
$${}^pR^mf_\ast {Q}_X = \sum_k \IC V^{m}_k[k]$$ 
\endproclaim

\definition{4.8 Remark} The uniqueness of the subobjects
${}^p \tau_{\leq m} Rf_\ast {Q}_X=$ \linebreak 
$\sum_{i\leq m} {}^pR^if_\ast {Q}_X [-i]$
is an immediate consequence of the axiom $\Ho (\De^{\leq 0}, \De^{>0})=0$
for $t$-structures. 
\enddefinition

\subhead {} Grothendieck motives over $S$ \endsubhead

Let $\underline M \in A \M(S)$ 
and $M =\r \underline M \in \De^b_{cc}(S)$ be its realisation.
Recall that $A \M (S)$ is the category of homological motives over $S$,
constructed in \S 2.  
As a consequence of the topological decomposition theorem,
$M$ is naturally equipped with an 
increasing filtration (perverse Leray filtration)
$$...\subset L_m M \subset L_{m+1} M \subset ... \subset M$$
defined as
$$L_m M = {}^p \tau_{\leq m} M$$ 
with $gr^L_m M = {}^p{\Cal H}^m (M)[-m]$. 
Let ${\underline u} : \underline M \to
{\underline N}$ be a morphism in $A \M(S)$ and $u=\r {\underline u}: M
\to N$ be its realisation in $\De^b_{cc}(S)$. Because $\r$ and ${}^p {\Cal
H}^m$ are functors, we get a system of compatible morphisms
$$L_m u : L_m M \to  L_m N$$
To elaborate more on this point, choose decompositions
$$\align 
M &=\sum {}^p{\Cal H}^m M[-m] \\ 
N &=\sum {}^p{\Cal H}^m N[-m]
\endalign$$
With respect to these decompositions, $u= \sum u^l{}_m$, where
$u^l{}_m: {}^p{\Cal H}^m M[-m]\to {}^p{\Cal H}^l N[-l]$. Here $u^l{}_m$ can
be regarded as an extension of perverse sheaves
$$u^l{}_m \in \Ho_{\De^b_{cc}S}\bigl({}^p{\Cal H}^m M[-m],
{}^p{\Cal H}^l N[-l] \bigr)=
\operatorname{Ext}^{m-l}_{\Pe S}\bigl({}^p{\Cal H}^m H,
{}^p{\Cal H}^l K \bigr)$$
hence $u^l{}_m = 0$ if $m<l$, corresponding to the fact that we have morphisms
$L_m u : L_m M \to  L_m N$,
and $u$ can be therefore visualised as an upper
triangular matrix
$$u=\pmatrix
   &\vdots &\vdots &   \\        
...&u^1{}_1&u^1{}_2&...\\
...&0      &u^2{}_2&...\\
...&0      &0      &...\\
   &\vdots &\vdots &
\endpmatrix$$

The above considerations imply that, passing to the corresponding graded
objects, we have a perverse realisation functor
$${}^p \r : A\M (S) \to gr\Pe (S)$$

\definition{4.9 Definition} The category $\M (S)$ of Grothendieck motives
over $S$ has the same objects as $A\M (S)$ and morphisms
$$\Ho_{\M S} ({\underline M}, {\underline N})=
\Img \bigl(\Ho_{A\M S} ({\underline M}, {\underline N})\to 
\Ho_{gr\Pe S} ({}^p\r{\underline M}, {}^p\r {\underline N})\bigr)$$
\enddefinition

\subhead {} Semisimplicity and decomposition in $\M (S)$ \endsubhead

In this subsection we assume that 
desingularisations of varieties over $k$ exist.

\definition{4.10 Notation}
 
(1) {\it In this subsection only}, underlined capital 
letters $\underline M$ denote
objects in $\M (S)$, while non underlined letters $M$ denote the corresponding
realisation in $gr \Pe (S)$. The same letter will denote a morphism
in $\M S$ or its realisation in $gr \Pe (S)$: the context will always
make it clear which is meant. When we want to specifically emphasise
the realisation functor, we call it ${}^p \r : \M (S) \to gr \Pe (S)$.

(2) If $X$ is a smooth variety and $f: X \to S$ a
projective morphism, we denote
$${}^p\underline R f_\ast Q_X$$
the corresponding object in $\M (S)$.
\enddefinition

We will prove the following results:

\proclaim{4.11 Theorem} Assuming the standard conjectures, 
the category $\M (S)$ is abelian and semisimple.
\endproclaim

\proclaim{4.12 Decomposition theorem in $\M (S)$} Assume
the standard conjectures.
Let $X$ be a smooth
variety and $f: X \to S$ be a projective morphism. 

(1) There is a canonical direct sum decomposition in $\M (S)$
$${}^p \underline R f_\ast {Q}_X \cong \sum {}^p\underline R^mf_\ast 
{Q}_X [-m]$$
where ${}^p\underline R^mf_\ast 
{Q}_X [-m]$ denotes an object in $\M (S)$, together with a given
isomorphism in $gr \Pe (S)$
$${}^p\r {}^p\underline R^mf_\ast {Q}_X [-m] @> \cong >> 
{}^p R^mf_\ast {Q}_X [-m]$$

(2) There is a canonical direct 
sum decomposition:
$${}^p \underline R^mf_\ast {Q}_X[-m] = \sum_k \underline \IC V^{m}_k[k-m]$$ 
where $V_k$ is a local system on a Zariski locally closed 
subvariety $T_k \subset S$ and 
$\underline \IC V^{m}_k[k-m]$ denotes an object in $\M (S)$, 
together with a given 
isomorphism in $gr \Pe (S)$
$${}^p\r\underline \IC V^{m}_k[k-m] @> \cong >> \IC V^{m}_k[k-m]$$
\endproclaim

\definition{4.13 Remark} The decomposition is unique (contrary to 4.7)
because of the way $\M (S)$ is built as a faithful subcategory 
of $grPerv (S)$.
\enddefinition

4.11 is an immediate consequence of the following (as is the case for
motives over a point, compare \cite{Kl1}) proposition 4.14 which 
will be shown, together with theorem 4.12, at the very end of this 
subsection.

\proclaim{4.14 Proposition} Let $X$ be a smooth variety and $f: X \to S$ a
projective morphism. Assuming the standard conjectures, 
$\En_{\M S } {}^p \underline R f_\ast Q_X$ is a semisimple ring, finite 
dimensional over ${\Bbb Q}$.
\endproclaim

The proof of the proposition will depend on the following:

\proclaim{4.15 Lemma (decomposition mechanism)} Let ${\Cal A}$ be an
abelian semisimple category, ${\Cal B}$ an additive category and
${\Cal A} \subset {\Cal B}$ a fully faithful embedding. Assume given
objects $A$, $A^\prime$ of ${\Cal A}$ and $B$ of ${\Cal B}$, and morphisms
$$A @> i_\ast >> B @> i^\ast >>A^\prime$$
denote $a=i^\ast i_\ast : A \to A^\prime$: by assumption this is a morphism
in ${\Cal A}$. 

There is a non unique projector $\beta : B \to B$, with image in ${\Cal
B}$, and a natural isomorphism
$$\varphi (\beta): \Img \beta @> \cong >> \Img a$$
\endproclaim

\demo{Proof} Let $V = \Img a$, $t: A \to V$ and $t^\prime : V \to A^\prime$
the natural maps. Choose $s:V \to A$, $s^\prime : A^\prime \to V$ such that
$$ts=Id_V, \quad s^\prime t^\prime = Id_V$$
and let $\alpha =s s^\prime : A^\prime \to A$. It is immediate to verify 
that $\alpha a \alpha =\alpha$ and $a \alpha a=a$. Let now $\beta =
i_\ast \alpha i^\ast$, it is immediate that $\beta^2 = \beta$.

{\smc Claim.} With the diagram
$$\diagram
B \rrto^\beta \drto_{s^\prime i^\ast} & & B \\
                 & V\urto_{i_\ast s}  & \enddiagram$$
$V= \Img \beta$. Let $X$ be any object of $\Cal B$. We wish to check that
$\Ho (X, B)$, resp. $\Ho (B, X)$ is an image of $\_\circ \beta$, resp.
$\beta \circ \_$, in the category of
{\it abelian groups}, via the diagram
$$\diagram
\Ho (X,B) \rrto \drto & &\Ho (X, B) \\
                 & \Ho (X, V) \urto  & \enddiagram$$
resp. the diagram
$$\diagram
\Ho (B, X) \rrto \drto & & \Ho (B, X) \\
                 & \Ho (V, X)\urto  & \enddiagram$$
In other words, {\it we may assume that} ${\Cal B}$ 
{\it is an abelian category}. The claim then follows from the observation that
$s^\prime i^\ast$ is surjective (indeed $s^\prime i^\ast i_\ast =t$ is 
surjective) and $i_\ast s$ is injective (indeed $i^\ast i_\ast s=t^\prime$
is injective).\qed \enddemo

\definition{4.16 Example (conic bundles)} As an example, we use the 
decomposition mechanism to {\it very briefly} outline the 
calculation of the Chow motive of a conic bundle, following \cite{Be}.

Let $X$ be a smooth 3-fold, $f: X \to S$ a conic bundle structure.
We assume, for simplicity, that the discriminant $\Delta \subset S$
is a smooth divisor, denote $Y=f^{-1} \Delta$ and $i:Y^\prime \to X$
the normalisation of $Y$. In an obvious way $Y^\prime \to \Delta$
factors through a ${\Bbb P}^1$-bundle $p: Y^\prime \to D$, with
a distinguished section (the conductor) $s: D \to Y^\prime$, by which
we may think $D \subset Y^\prime$, and
an \'etale double cover $D \to \Delta$. Let $\tau : D \to D$ be the
involution associated to this double cover. 

{\smc Step 1.} Let $$a= i^\ast i_\ast : h_S Y^\prime (-1) 
\to h_S Y^\prime $$ in $CH \M(S)$, then
$$a=c_1(N_\Delta S)-s_\ast(s^\ast-\tau_\ast s^\ast)$$
where, abusing notation slightly, $N_\Delta S$ is the pull back to
$Y$ of the normal bundle $N_\Delta S$ of $\Delta$ in $S$.
Indeed, let $\Gamma = \Gamma_i\subset Y^\prime \times X$ be the graph of
$i$, and $\gamma \in CH_2(Y^\prime \times X)$ its class. Then $i_\ast = \gamma$
and $i^\ast ={}^t\gamma$, while by definition $a=p_{13 \ast}
(p_{23}^\ast {}^t \gamma \cdot p_{12}^\ast \gamma)$ 
(definition 2.1). We can calculate the 
intersection product $p_{23}^\ast {}^t \gamma \cdot p_{12}^\ast \gamma$
with the help of the following fibre square
diagram
$$\diagram 
D \coprod  Y^\prime \rrto \dto_{(s, s\tau)\coprod \Delta_{Y^\prime}} & 
&Y^\prime \times Y^\prime \dto^{(1,i)\times 1} \\
Y^\prime \times Y^\prime \rrto^{1\times (i,1)} & &Y^\prime \times X \times 
Y^\prime\enddiagram$$
By the excess intersection formula then
$$a=c_1(E)+s_\ast\tau_\ast s^\ast$$
where $E$ is the excess bundle on $Y^\prime$, defined by the exact sequence
$$0\to T_{Y^\prime}\to i^\ast T_X \to E \to 0$$
Finally, it is easy to convince oneself that $E=N_\Delta S(-D)$, 
giving $c_1(E)=$ \linebreak 
$c_1(N_\Delta S)-D=c_1(N_\Delta S)-s_\ast s^\ast$, which proves
our formula. 

{\smc Step 2.} Now $Y^\prime \to D$ is a ${\Bbb P}^1$-bundle, therefore
we have an isomorphism
$$(p^\ast , s_\ast ): 
h_S D \oplus h_S D (-1)\to
h_S Y^\prime$$
Using this isomorphism, it is quite easy to see that
$$a: h_S D(-1) \oplus h_S D(-2) \to
h_S D \oplus h_S D(-1) $$
can be written in matrix form as
$$a=\pmatrix
c_1(N_\Delta S) & 1-\tau^\ast \\
0 & c_1(N_\Delta S \otimes N_D^\vee Y^\prime \otimes 
\tau^\ast N_D Y^\prime)\endpmatrix$$

{\smc Step 3.} From the previous step and the decomposition
mechanism we can see, for instance, the classical result 
stating that
the Prym motive $(h^1 D, 1-\tau^\ast)$ is a direct summand of 
the intermediate motive $h^3 X (1)$. 
\enddefinition

Before we embark in the proof of proposition 4.14 and Theorem 4.12, we need 
a definition and a lemma.

\definition{4.17 Definition} Let $X$ be a smooth variety, $f: X \to S$ a
projective morphism. An {\it equisingular stratification} of $f$ is a 
pair $({\Cal T},{\Cal Y}^1 \rightrightarrows {\Cal Y}^0 \to {\Cal Y})$ where:

(1) ${\Cal T}=\{T_k\}$ is a good stratification of $S$,

(2) ${\Cal Y}=\{Y_k\}$ and $Y_k$ is defined by the fibre square
$$\diagram
Y_k \dto  \rto &X \dto \\
T_k\rto           &S\enddiagram$$

(3) ${\Cal Y}^1=\{Y_k^1\}$, ${\Cal Y}^0=\{Y_k^0\}$ and $Y_k^1
\rightrightarrows Y_k^0 \rightarrow Y_k$ is a truncated simplicial resolution.

The above data are subjected to the following condition:

(4) The compositions $Y_k^1 \to T_k$, $Y_k^0 \to T_k$ are all smooth
(not necessarily equidimensional). In particular, for all $t \in T_k$
$$Y^1_{k,t} \rightrightarrows Y^0_{k,t} \rightarrow Y_{k,t}$$
is a truncated simplicial resolution.
\enddefinition

It is a consequence of our assumption on the existence of resolutions
of singularities, that
equisingular stratifications of $f: X \to S$ exist. We will need the
following:

\proclaim{4.18 Key lemma} 
Let $X$ be a smooth variety, $f: X \to S$ a projective
morphism. Fix an equisingular stratification $({\Cal T},
{\Cal Y}_1 \rightrightarrows {\Cal Y}_0 \to {\Cal Y})$ of the morphism
$f$. Let $T_0$ be the smallest stratum of ${\Cal T}$. We have a natural 
isomorphism (cf. the notation in the statement of
the topological decomposition theorem 4.7)
$$i_{T_0 \ast} V^{m}_0= 
\Img \bigl( R^{m} f_\ast i_{Y_0 \ast} i_{Y_0}^! {Q_X} 
\to R^{m}f_\ast {Q}_X\bigr)$$
\endproclaim

\demo{Proof} {\smc Step 1.} First of all, if $h<k$ the natural maps
$${\Cal H}^j i_{{T_h}\ast} i_{T_h}^! IC V^{m}_k \to {\Cal H}^j 
IC V^{m}_k$$
are zero for all $j$. 

{\bf Warning:} this is not saying that
$i_{{T_h}\ast} i_{T_h}^! \IC V^{m}_k \to 
\IC V^{m}_k$ is the zero map in $\De^b_{cc}(S)$. 

Indeed:
${\Cal H}^j i_{T_h}^! \IC V^{m}_k=0$ for $j\leq k-h$, and
${\Cal H}^j i_{T_h}^\ast \IC V^{m}_k=0$ for $j\geq k-h$.

{\smc Step 2.} We have a fibre square
$$\diagram
Y_0 \rto^{i_{Y_0}} \dto & X \dto^f \\
T_0 \rto^{i_{T_0}} &S \enddiagram$$ 
From it, using the proper base change theorem and the
topological decomposition theorem, we derive the following commutative diagram
$$\diagram
Rf_\ast i_{Y_0\ast} i_{Y_0}^! {Q}_X \ddouble \rto 
& Rf_\ast {Q}_X \ddouble \\
i_{T_0\ast}i_{T_0}^!Rf_\ast{Q}_X \ddouble\rto 
& Rf_\ast {Q}_X \ddouble \\
\overset \sum_{k>0} \sum_m i_{T_0\ast} i_{T_0}^!\IC V^{m}_k[k-m]
\to{\underset 
\sum_m V^{m}_0 [-m] \to
\oplus} \rto &
\overset \sum_{k>0} \sum_m \IC V^{m}_k [k-m]\to
{\underset \sum_m V^{m}_0 [-m]   \to
\oplus} 
\enddiagram$$
The result then follows from step 1, upon taking ${\Cal H}^{m}$ of both sides
of the bottom portion of the diagram. \qed \enddemo

\subsubhead Proof of proposition 4.14 and theorem 4.12 \endsubsubhead
 
Fix an equisingular stratification $({\Cal T},
{\Cal Y}^1 \rightrightarrows {\Cal Y}^0 \to {\Cal Y})$ of the morphism $f:X
\to S$. The proof is by induction on $\di S$ and the number of strata in
${\Cal T}$. The basis for the induction is solid because:

(a) If $\di S=0$, $\M (S) = \M$ is semisimple by 3.5, proven in \cite{Kl1}.
The decomposition theorem 4.12 in this case can be proven as follows.
Again in \cite{Kl1} is shown that, if the standard conjecture of Lefschetz
type holds, then there are cycles $\Pi^i$ representing the K\"unneth
components $\pi^i$ of the diagonal $\Delta \subset X \times X$, and the sought
for decomposition is then
$$X = \sum (X, \Pi^i)$$

(b) If ${\Cal T}$ has only one stratum, note that, by definition
of equisingular stratification, this happens if and only if 
$f: X \to S$ is a smooth morphism. By 3.4, there is a cycle $Z$ on
$X \times_S X$ inducing the $\Lambda$ operator on each fibre.
The same proof as in \cite{Kl1} will then show that $\En_{\M S} X$
is a semisimple ring, finite dimensional over ${\Bbb Q}$. Then
again, as in \cite{Kl1}, there are classes $\Pi^i \in CH_{\di X}
X \times_S X$ inducing on each fibre the K\"unneth components of the diagonal
of that fibre, and one can get the decomposition as above
$$X = \sum (X, \Pi^i)$$

Let now $T_0$ be the smallest stratum. We apply the decomposition mechanism
4.15 
with the following setup. ${\Cal A}=\M (T_0)$, which is by inductive
assumption abelian and semisimple, ${\Cal B}=\M (S)$. 
We take
$$\align
\underline A&= \Cok \bigl({}^p\underline Rf_\ast D_{Y^1_0}(-\di X) \to 
{}^p\underline Rf_\ast D_{Y^0_0} (-\di X) \bigr) \\
\underline A^\prime&=\Ker\bigl( {}^p\underline  Rf_\ast Q_{Y^0_0} \to {}^p
\underline Rf_\ast Q_{Y^1_0}\bigr) 
\endalign$$
and
$$\underline B={}^p\underline Rf_\ast Q_X$$
$\underline A$ and $\underline A^\prime$ are objects of $\M (T_0)$, but 
if we like we can think of them as being in $\M(S)$ via the obvious
inclusion $\M (T_0) \subset \M (S)$.
There are obvious maps $i_\ast : \underline A \to \underline B$ and 
$i^\ast : \underline B \to \underline A^\prime$. We
will not need this, but we still like to say that the assignment $X$ to 
$\underline A$,
resp. $\underline A^\prime$ is functorial, in other words it does not depend 
on the choice of the equisingular stratification.
As we anticipated, we will denote $A$, $A^\prime$ and $B$ the realisations
in $gr \Pe  (S)$.

{\smc Claim 1.} Let $R^{m}={\Cal H}^{m} A$, $R^{m, \prime}={\Cal
H}^{m} A^\prime$. Then there is a natural isomorphism
$$V^{m}_0=
\Img \bigl( {\Cal H}^{m} a: R^{m} \to R^{m, \prime}\bigr)$$
In fact we know from the key lemma that
$$V^{m}_0=
\Img \bigl( {\Cal H}^{m}  i_{T_0}^! Rf_\ast {Q}_X \to
{\Cal H}^{m} i_{T_0}^\ast Rf_\ast{Q}_X \bigr)$$
Since ${\Cal H}^{m} i_{T_0}^! Rf_\ast {Q}_X$ (resp. 
${\Cal H}^{m} i_{T_0}^\ast Rf_\ast {Q}_X$) has weights $\geq m$
(resp. $\leq m$), we have
$$i_{T_0 \ast} V^{m}_0=
\Img \bigl( gr^W_{m} {\Cal H}^{m} i_{T_0}^! Rf_\ast {Q}_X \to
gr^W_{m} {\Cal H}^{m} i_{T_0}^\ast Rf_\ast {Q}_X \bigr)$$ 
The claim now follows from the identifications
$$\diagram
gr^W_{m}{\Cal H}^{m} i_{T_0}^! Rf_\ast {Q}_X \ddouble \rto &
gr^W_{m}{\Cal H}^{m} i_{T_0}^\ast Rf_\ast {Q}_X \ddouble\\
R^{m} \rto^{{\Cal H}^{m} a} &R^{m, \prime} \enddiagram$$

The decomposition mechanism, together with the claim, provides a projector 
$\beta \in \En_{\M S } {}^p\underline Rf_\ast Q_X$ s.t., upon setting
${\underline M} = \Ker \beta$, ${\underline V}=\Img \beta$, we have
$${}^p\underline Rf_\ast Q_X = {\underline M} \oplus {\underline V}$$
whose realisations in $gr \Pe (S)$ are
$$\align
V&=\sum_m i_{T_0 \ast} V^{m}_0[-m]\\
M&=\sum_{k>0} \sum_m \IC V^{m}_k[k-m]\endalign$$

Let us now prove 4.14, i.e.,
let us show that $\En_{\M S} {}^p\underline Rf_\ast Q_X$ is a 
semisimple ring,
finite dimensional over ${\Bbb Q}$. According to
the above decomposition
$$\En_{\M S} {}^p \underline R f_\ast Q_X =\En_{\M S} 
{\underline M} \oplus
\En_{\M T_0} {\underline V}$$
Now $\En_{\M T_0} {\underline V}$ is semisimple finite dimensional over 
${\Bbb Q}$ by inductive assumption on $\dim S$.
The same is true of $\En_{\M S} {\underline M}$ by inductive
assumption on the number of strata, by the 

{\smc Claim 2.}
$$\En_{\M S} {\underline M}=
\En_{\M (S\setminus T_0)} ({\underline M}|_{S\setminus T_0})$$

Indeed we have a diagram
$$\diagram
\En_{\M S} \bigl({\underline M}\bigr) \dto^{\text{surj}} \rrto^{\text{inj}}&
&\En_{gr\Pe S}\bigl(\sum_{m, k>0} \IC V^{m}_k\bigr)\ddouble\\
\En_{\M (S\setminus T_0)} \bigl(
{\underline M}|_{S\setminus T_0}\bigr)\rrto^{\text{inj}}& &
\En_{gr\Pe (S\setminus T_0)}\bigl(\sum_{m, k>0} 
\IC V^{m}_k|_{S\setminus T_0}\bigr)\enddiagram$$
where surjectivity of the left vertical arrow follows from 
the fact that cycles can be closed:
$$\diagram
CH_\bullet X \times_S X \rto \dto^{\text{surj}}&H^{BM}_{2\bullet} X\times_S X
\dto\\
CH_\bullet X\times_{S\setminus T_0}X\rto&H^{BM}_{2\bullet}
X\times_{S\setminus T_0}X\enddiagram$$

This finishes the proof of 4.14, beacuse $\En_{\M S} X$ is direct sum of 2 
rings, each of which is semisimple and finite dimensional over ${\Bbb Q}$.

Finally, we shall now prove the decomposition theorem 4.12. By induction on
the number of strata, we may assume that the decomposition theorem holds
over $S\setminus T_0$: 
$${}^p\underline R f_\ast Q_X |_{S \setminus T_0}=
\sum_{k>0} \sum_m \underline \IC V_k^m[k-m]$$
To give such a decomposition is equivalent to giving the projectors
$\Pi_k^m \in$ \linebreak
$\En_{\M (S\setminus T_0)}$ down to $\underline \IC V_k^m$.
Since by claim 2 $\En_{\M S} {\underline M}=
\En_{\M (S\setminus T_0)} ({\underline M}|_{S\setminus T_0})$, this
also decomposes ${}^p\underline R f_\ast Q_X$ over $S$ into the desired
pieces.

\head 5. Filtrations for quasiprojective varieties \endhead
 
Throughout this section we assume that
desingularisations of varieties over $k$ exist.  
For a quasiprojective variety $X$ we will put a canonical filtration on 
its rational Chow group $CH_s X$ so that some functorial properties are 
satisfied, theorem 5.1. The filtration is shown 
to satisfy additional properties
if one assumes the conjectures of Grothendieck and Murre, 5.2.  

As an application, assuming these conjectures, we show that the 
projectors in \S 4 can be lifted to an orthogonal set of projectors in
the Chow group of 
relative self correspondences, 5.10. 
This is done by showing that the map 
$$CH_{\dim X} X\times_S X \to \En_{\M S} {}^p \underline Rf_* {Q}_X$$
is a surjective homomorphism with nilpotent kernel. 
\smallskip 

\proclaim{5.1 Theorem}
For a quasi-projective variety $X$, there is a decreasing finite filtration 
{\rm (the canonical filtration)} on its Chow group $CH_s X$  
$$CH_s X=F^0 CH_s X\supset F^1 CH_s X \supset F^2 CH_s X\supset 
\cdots $$
subject to the following conditions: 

(i) If $f: X\to Y$  is a proper map of quasi-projective varieties, then the induced map
 $f_* : CH_s X\to CH_s Y$ respects the filtrations, i.e. 
$$f_* F^\nu CH_s X\subset F^\nu CH_s Y$$
for each $\nu$.  If $f$ is proper and surjective, then $f_* F^\nu CH_s X = F^\nu CH_s Y$
(in other words, the surjection $f_*$ is {\it strictly} compatible with $F^\bullet$).  

(ii)  If $j : U\to X$ is an open immersion of quasi-projective varieties, then the restriction $j^*: CH_s X\to 
CH_s U$ is strictly compatible with $F^\bullet$: $j^*F^\nu  CH_s X = F^\nu  CH_s U$. 

(iii) For a smooth projective $X$
$$F^1 CH_s X = CH_s (X)_{hom}= \Ker \bigl(cl: CH_s X\to H _{2s} X\bigr)$$ 
where $cl$ is the cycle class map. 

(iv)  The external product map 
$CH_s X\otimes CH_t Y\to CH_{s+t} (X\times Y)$ respects $F^\bullet$: if $z\in F^\nu CH_s X$ and 
$w\in F^\mu  CH_t Y$ then $z\times w\in F^{\nu+\mu} CH_{s+t} (X\times Y)$. 

(v)  The internal product respects $F^\bullet$: if $X$ smooth quasiprojective equidimensional, $z\in 
 F^\nu CH_s X$ and 
$w\in F^\mu  CH_t X$ then $z \cdot  w\in F^{\nu+\mu} CH_{s+t -\dim X} X$. 

(vi) Refined Gysin maps respect $F^\bullet$.  Let $i : X\to Y$ be a 
regular embedding of codimension $d$, and 
$$\diagram
X^\prime \rto \dto &Y^\prime \dto \\
X\rto_i & Y\enddiagram$$
be a Cartesian square where $Y^\prime$ is an arbitrary quasiprojective variety 
and $Y^\prime \to Y$ is an arbitrary map. Then the refined Gysin map \cite{Fu, 
Chap. 6}
$$i^! : CH_s Y'\to CH_{s-d} X'$$
respects $F^\bullet$. 

(vii)  If $X, Y$ are quasi-projective varieties, $X$ equidimensional, and $p_Y: X\times Y\to Y$ is the projection, then the map 
$p_Y^*: CH_s Y\to CH_{s+\dim X} (X\times Y)$ respects $F^\bullet$. 

(viii) Let $W, X$ be smooth projective equidimensional,  
$\Gamma \in CH_{\dim X-i} (W\times X)_{hom}$, and 
$$\Gamma_* : CH_{s+i} W\to CH_s X$$
the induced map (by (i), (v) and (vii), $\Gamma_* $ respects $F^\bullet$). 
Then the map induces zero on the $F$-graded pieces:
$$Gr_F^\bullet \Gamma_* =0 : 
Gr_F^\bullet  CH_{s+i} W\to Gr_F^\bullet CH_s X$$
\endproclaim 

\definition{5.2 Remarks}

(1)  As already noted in 3.14 and 3.18, S. Saito defined a filtration 
 $F^\bullet$ on $CH_s (X)$ for $X$ smooth projective satisfying 
the conditions (i), (iii), (v), (vii) (where $X, Y$ are smooth projective) and 
(viii). 

If we further assume Murre's conjecture, it follows that Saito's 
filtration is separated, i.e. for any $X$ and $s$, one has 
$F^\nu CH_s X=0$ for $\nu$ large.
 
In the proof of Theorem 5.1, we take Saito's filtration and show that 
it uniquely extends to a filtration for $X$ quasiprojective, so that the 
conditions (i)-(viii) are satisfied. 

(2) (iv) and (v) follow from (vi) and (vii). 

(3) We will use the following properties of refined Gysin maps 
\cite{Fu, Chap. 6}:
compatibility with proper push forwards, compatibility with flat pull-backs, 
the excess intersection formula, and the fact: if $i$ is of codimension one, 
namely if $X\subset Y$ is a Cartier divisor, then $i^!$ coincides with 
intersection with $X$. 

(4) There is an interpretation of the filtrations in terms of 
mixed motives; we will not need this. 
\enddefinition

\proclaim{5.3 Theorem} Assuming Grothendieck's and Murre's conjectures,
the filtration in theorem 5.1 satisfies in addition the following properties:

(ix)  For any quasiprojective variety 
$$F^1 CH_s X = CH_s (X)_{hom}= \Ker \, \bigl(cl: CH_s X\to H^{BM} _{2s} X
\bigr)$$
where $cl: CH_s X\to H^{BM} _{2s} X$ is the cycle map into Borel-Moore 
homology. 

(x) For each $X$, one has  $F^\nu CH_s X= 0$ for $\nu$ large. 
\endproclaim

To show Theorem 5.1, we 
take Saito's filtration $CH_s X$ for $X$ 
smooth projective and attempt to extend it to $X$ general.  

First consider the case $X$ is smooth quasiprojective. Take a smooth 
projective variety $\overline X$ 
and an open immersion $j : X\hookrightarrow \overline X$.  
It induces the surjective map 
$$j^* : CH_s \overline X\to CH_s X$$
and  $CH_s X$ is given the induced filtration: $F^\nu CH_s X = j^* F^\nu 
CH_s \overline X $. 
This filtration is independent of the choice of a compactification. 
In fact, let 
$j' : X\to \overline{X}^\prime$ be another smooth compactification. Since 
$\overline X$ and $\overline{X}^\prime$ are dominated 
by a third compactification, one may assume that there is a map $f: 
\overline{X}^\prime\to \overline X$ such that 
$f\scirc j'= j$. The diagram 
$$\diagram
CH_s \overline{X}^\prime \ddto_{f_\ast} \drto^{j^{\prime \ast}}& \\
         &CH_s X \\
CH_s \overline X \urto_{j^\ast} & \enddiagram$$
commutes.  By the strictness of $f_*$   (5.1 (i) for surjective maps
of smooth projective varieties) one has 
$$f_* F^\nu CH_s \overline{X}^\prime = F^\nu CH_s \overline X$$
so $j^*$ and $j'{}^*$ induce the same filtrations. 
\smallskip 

\proclaim{5.4 Proposition} 

(1) If $X$ and $Y$ are smooth quasiprojective varieties and 
$f: X\to Y$ is proper (resp. proper surjective) then 
the map $f_*: CH_s X\to CH_s Y$ respects $F^\bullet$ (resp. strictly 
compatible with $F^\bullet$). 

(2) If $j: U\hookrightarrow X$ is an open immersion of 
smooth varieties, $j^* : CH_s X
\to CH_s U$ is strictly compatible with $F^\bullet$. 
\endproclaim

\demo{Proof} 
For (1) we take smooth compactifications $\overline X$, $\overline Y$ of $X$, $Y$, respectively, so that $f$ extends to a map $\overline f : \overline X\to \overline Y$. 

Consider the commutative diagram 
$$\diagram
CH_s X \rto^{f_\ast} & CH_s Y \\
CH_s \overline X \uto \rto_{\overline f_\ast} &CH_s \overline Y 
\uto \enddiagram$$
where the vertical arrows are the pull backs by open immersions. Since $\overline f_*$ respects $F^\bullet$ (resp. strictly compatible with $F^\bullet$ if $f$ is surjective)
$\overline f_* F^\nu CH_s \overline X \subset  F^\nu CH_s \overline Y$ 
(resp. equal).  On the other hand, by definition $F^\nu CH_s \overline X$ surjects to 
$F^\nu CH_s X$, and $F^\nu CH_s \overline Y$ surjects to 
$F^\nu CH_s Y$. 
 Hence $ f_* F^\nu CH_s X \subset F^\nu CH_s Y$ (resp. equal). 

For (2) take a compactification 
$j: X\to \overline X$ and consider the commutative diagram 
$$\diagram
CH_s X\rrto^{j^\ast} &  &CH_s U                    \\
            &CH_s \overline X \ulto \urto \enddiagram$$
where all the arrow are restrictions by open immersions.  Since $F^\nu CH_s 
\overline X$ surjects to both 
$F^\nu CH_s X$ and  $F^\nu CH_s U$, one has $j^* F^\nu CH_s X = F^\nu CH_s U$. \qed \enddemo

For an arbitrary quasiprojective variety $X$, take a desingularization $\pi: \tilde X\to X$ and equip 
$CH_s X$ with the filtration induced by the surjective map $\pi _* : CH_s 
\tilde X\to CH_s X$. 
By an argument using 5.4(1), one sees that the filtration is 
well defined independent of the choice of $\tilde X$.

\proclaim{5.5 Proposition}

(1) If $X$, $Y$ are quasi-projective varieties and $f: X\to Y$ is proper 
(resp. proper surjective), then $f_* : CH_s X\to 
CH_s Y$ respects $F^\bullet$ (resp. strictly compatible with $F^\bullet$). 

(2) If $j: U\to X$ is an open immersion of quasi-projective varieties, 
$j^* : CH_s X\to CH_s U$ is strictly compatible with $F^\bullet$. 
\endproclaim

\demo{Proof} To prove (1), take desingularizations 
$\pi:\tilde X\to X$, $\pi':\tilde Y\to Y$ so that $f$ extends
to a map $\tilde f : \tilde X\to \tilde Y$.  Then 
$\pi'\scirc \tilde f= f\scirc \pi$ and one has a 
commutative diagram 
$$\diagram
CH_s X \rto^{f_\ast} & CH_s Y \\
CH_s \tilde X \uto^{\pi_\ast} \rto_{\tilde f_\ast} &CH_s\tilde Y 
\uto_{\pi^\prime_\ast}\enddiagram$$
where the arrows are proper push forwards. Since the map 
$\tilde f_*$ respects $F^\bullet$ (resp. stricly compatible 
with $F^\bullet$) by 5.4(1), and so are the vertical surjective maps by definition, $f_*$ respects $F^\bullet$ (resp. is strictly compatible with $F^\bullet$).  
The proof of (2) is similar. \qed \enddemo

\demo{Proof of Theorem 5.1}
The properties (i) and (ii) have been verified, and we 
started with the filtration satisfying (iii) and (viii). 

(vii) This is verified by reducing first to the case where $X$ and $Y$ are 
both smooth, and then to the case they are smooth projective. 

(vi) Follows from strictness of the filtration under proper
maps and open immersions (i) and (ii), lemma 5.6 below and compatibility
of the filtration under action of correspondences
on smooth projective varieties (i), (v), (vii). \qed \enddemo 

\proclaim{5.6 Lemma} Let $i:T \to S$ be a regular embedding of codimension 
$d$
$$\diagram
Y \dto_g \rto & X \dto^f \\
T \rto_i      & S \enddiagram$$
a fibre square, and
$$i^!: CH_r X \to CH_{r-d}Y$$
the associated refined Gysin map \cite{Fu, Ch. 6}. There are: 

(1) smooth varieties
$U$, $V$, proper surjective maps $p: U \to X$, $q:V \to Y$, 

(2) a commutative diagram
$$\diagram
V \dto_q \rto & U \dto^p \\
Y \rto      & X \enddiagram$$

(3) smooth compactifications $U \subset \overline U$, $V \subset \overline V$
and a correspondence $\Gamma \in CH_\bullet \overline V \times \overline U$,
such that for a cycle $\alpha \in CH_r \overline U$
$$i^! p_\ast \bigl(\alpha|U\bigr) 
=q_\ast \bigl((\Gamma^\ast \alpha)|V\bigr)$$
\endproclaim

\demo{Proof} 

{\smc Step 1.} In this step we reduce the problem to the case
where $Y=Y_1\coprod Y_2\hookrightarrow X$, $X$ is smooth, 
$Y_1 \hookrightarrow X$ is a normal crossing divisor, and 
$Y_2 \hookrightarrow X$ is the inclusion of a bunch of connected components.

Let indeed $\delta: X^\prime \to X$ be a resolution of singularities such that
$Y^\prime = Y \times_X X^\prime =Y^\prime_1 \coprod Y^\prime_2
\hookrightarrow X^\prime$ is as above. 
We have a commutative diagram of fibre squares
$$\diagram
Y^\prime \dto_\varepsilon \rto & X^\prime \dto^\delta \\
Y \dto_g \rto                  & X \dto^f             \\
T \rto_i                       & S \enddiagram$$
By compatibility of refined Gysin maps with proper push forward \cite{Fu, Ch 6}
we have
$$i^!\delta_\ast = \varepsilon_\ast i^!$$
Therefore, it is enough to prove the result for $Y^\prime \to X^\prime$.

{\smc Step 2.} We now assume that
$j:Y\hookrightarrow X$ is a normal crossing divisor inside a smooth 
quasiprojective variety.
Let $E=g^\ast N_T S /N_Y X$ be the excess bundle
(by assumption, it has rank $d-1$). The excess intersection 
formula reads
$$i^! \alpha = c_{d-1} E \cap j^! \alpha$$
If now $U= X$, $U \subset \overline U$ is a smooth projective 
compactification of $U$ such that
the closure $\overline Y$ of $Y$ in $\overline U$ is a normal crossing
divisor,
then denoting $c_{d-1} \overline E$ any extension of $c_{d-1} E$
to $\overline U$,
$\nu: \overline V \to \overline Y$ the normalization, $V = \nu^{-1}
Y$ the normalization of $Y$,
and $h: \overline V \to \overline U$ the composition $\overline V \to 
\overline Y \to \overline U$, we have for a class
$\alpha \in CH_r \overline U$
$$j^! \bigl(\alpha|U\bigr)= \bigl(\nu|V\bigr)_\ast 
\bigl((h^\ast \alpha)|V\bigr)$$
and
$$i^! \bigl(\alpha|U\bigr) = 
\bigl(\nu|V\bigr)_\ast \bigl((c_{d-1} \overline E \cap h^\ast \alpha)|V
\bigr)$$
To conclude the proof now just take a correspondence
$$\Gamma \in CH_\bullet \overline V \times \overline U$$
such that
$$c_{d-1} \overline E \cap h^\ast \_     = \Gamma^\ast \_$$  

{\smc Step 3.} Finally, we treat the case where $j:Y\hookrightarrow X$ is
the inclusion of a bunch of connected components. This is quite a bit easier
than step 2: here
$$i^! \alpha =c_d (E) \cup j^\ast \alpha$$
We let $V=Y$, $U=X$,
$U \subset \overline U$ a smooth compactification and $c_{d} \overline E$ 
any extension of $c_{d} E$ to $\overline U$, 
$\overline j : \overline V \to \overline
U$ the corresponding compactification of $V \hookrightarrow U$. 
$\Gamma$ works if $\Gamma_\ast\_ =c_d \overline E \cup \overline j^\ast \_$.
\qed \enddemo

The following proposition will be needed in the proof of Theorem 5.3. 

\proclaim{5.7 Proposition} (Assume Grothendieck's standard conjectures.)
For a quasiprojective variety $X$, let 
$$H^{BM}_{2s} (X)_{alg}= \Img \bigl(cl: CH_s X\to H^{BM}_{2s} X\bigr)$$
which is a ${\Bbb Q}$-vector space. 

(1) If $j: U\hookrightarrow X$ is an open immersion of quasiprojective 
varieties and $i: Z= X-U\hookrightarrow X$ 
is the closed immersion of the complement 
of $U$, the exact sequence 
$$H^{BM}_{2s} Z\overset{i_*}\to\to H^{BM}_{2s} X\overset{j^*}\to\to H^{BM}_{2s} U$$
induces the following exact sequence on algebraic parts:
$$H^{BM}_{2s} (Z)_{alg}\overset{i_*}\to\to H^{BM}_{2s} (X)_{alg} \overset{j^*}\to\to H^{BM}_{2s} (U)_{alg}\to 0$$

(2) Let 
$$\diagram
Z^\prime \rto^{i^\prime} \dto_{q} & X^\prime \dto^p \\
Z \rto_i &X\enddiagram$$
be a Cartesian square of quasiprojective varieties such that the horizontal maps are closed immersions,
the vertical 
maps are proper surjective and $p$ induces an isomorphism 
$X'-Z'@>\cong >> X-Z$. 
Then the exact sequence 
$$H^{BM}_{2s} Z' @>i^\prime_\ast, q_\ast >> H^{BM}_{2s} X'\oplus  
H^{BM}_{2s} Z @>p_\ast -i_\ast >> H^{BM}_{2s} X$$
induces an exact sequence 
$$H^{BM}_{2s} (Z')_{alg} @>i^\prime_\ast, q_\ast >>
H^{BM}_{2s} (X')_{alg}\oplus  H^{BM}_{2s} (Z)_{alg}
@>p_\ast -i_\ast >>
H^{BM}_{2s} (X)_{alg}\to 0$$
\endproclaim

\demo{Proof} (1) 
We recall that $H^{BM}_i X$ for $X$ quasiprojective has a weight 
filtration $W_\bullet$, the weights are $\ge -i$, and the 
maps $i_*$, $j^*$ are strictly compatible with the weight filtrations. 

The first exact sequence induces, upon taking $Gr^W_{-2s}$, the exact sequence 
$$W_{-2s} H^{BM}_{2s} Z\to W_{-2s} H^{BM}_{2s} X \to    W_{-2s} H^{BM}_{2s} U$$
This may be viewed as an exact sequence in the category of Grothendieck 
motives $\M$. 
More specifically in the weight spectral sequence 
$${}_W E_1^{p\, q}\Rightarrow H^{BM}_{-p-q} X$$
which induces the weight filtration of $H^{BM}_{-p-q} X$, 
each ${}_W E_1^{p\, q}$ is the cohomology of a smooth projective variety and can be regarded as a Grothendieck motive denoted ${}_W \underline{E}_1^{p\, q}$; the differentials $d_1^{p, q}$ are morphisms of Grothendieck motives. Define 
$$Gr^W_{q} \underline{H}^{BM}_{-p-q} X= {}_W{\underline Z}_1^{p,q}
/{}_W{\underline B}_1^{p,q}$$
Taking the cohomological realisation $H^\ast$ we have 
$$H^*(Gr^W_{q} \underline{H}^{BM}_{-p-q} X) =Gr^W_{q} {H}^{BM}_{-p-q} X$$
because $H^\ast$ is exact (by the standard conjectures, $\M$ is a 
{\it semisimple} abelian category) and 
the weight spectral sequence above degenerates at $E_2$. 

Since $\M$ is semisimple, the functor 
$\Hom (pt(-s),   -): \M\to Vec_{{Q}}$ is also exact,  
and the above exact sequence induces an exact sequence
under this functor.  It is the desired exact sequence, except 
possibly at the end, because:
for a quasiprojective variety $X$
$$\Hom (pt(-s) , W_{-2s} \underline{H}^{BM}_{2s} X)=  H^{BM}_{2s} (X)_{alg}
$$
The surjectivity at the end is obvious since $j^*: CH_s X\to CH_s U$ is surjective. 

(2) One has  a commutative diagram 
$$\diagram
H^{BM}_{2s}Z^\prime\rto^{i^\prime_\ast}\dto^{q_\ast}&
H^{BM}_{2s}X^\prime\rto\dto^{p_\ast}&
H^{BM}_{2s}(X^\prime\setminus Z^\prime)\dto^{p_\ast}\\
H^{BM}_{2s}Z\rto_{i_\ast}&
H^{BM}_{2s}X\rto&
H^{BM}_{2s}(X\setminus Z)\enddiagram$$
where the rows are exact and the third vertical map is an isomorphism. 
Hence the exact sequence of Borel-Moore homology. 

Applying (1) to the closed immersions $i^\prime$ and $i$ respectively
one obtains a similar commutative diagram with exact rows
$$\diagram
H^{BM}_{2s}(Z^\prime)_{alg}\rto^{i^\prime_\ast}\dto^{q_\ast}&
H^{BM}_{2s}(X^\prime)_{alg}\rto\dto^{p_\ast}&
H^{BM}_{2s}(X^\prime\setminus Z^\prime)_{alg}\dto^{p_\ast}\\
H^{BM}_{2s}(Z)_{alg}\rto_{i_\ast}&
H^{BM}_{2s}(X)_{alg}\rto&
H^{BM}_{2s}(X\setminus Z)_{alg}\enddiagram$$
(the third vertical map is an isomorphism); the last exact sequence follows from this. \qed \enddemo 

\demo{Proof of Theorem 5.3}  
(x) is verified by reducing to the smooth projective case where it holds true, see 5.2(1). 

To show (ix)
first assume $X$ is smooth quasiprojective.  Take a compactification $j: X\to \overline X$, let 
$Z= \overline X-X$, and consider the commutative diagram with exact rows
5.7(1)
$$\diagram
H^{MB}_{2s}(Z)_{alg}\rto^{i_\ast}&
H^{MB}_{2s}(\overline X)_{alg}\rto^{j^\ast}&
H^{MB}_{2s}(X)_{alg}\rto &0 \\
CH_sZ\uto\rto_{i_\ast}&CH_s\overline X\uto\rto_{j^\ast}&CH_sX\uto\rto
&0\enddiagram$$
where the vertical arrows are the cycle class maps. 
Since $F^1 CH_s \overline X$ surjects to $F^1 CH_s X$ the 
 inclusion $F^1 CH_s X\subset  CH_s (X)_{hom}$ is obvious. The other inclusion follows from the above 
diagram. 

The proof in the general case is reduced to the smooth case by taking 
desingularisations and using 5.7(2). \qed \enddemo

As an application, let $S$ be a quasiprojective variety, $X$ a smooth 
quasiprojective equidimensional variety, and $f: X\to S$ a projective map. 
We recall that the group $CH_{\dim X} X\times _S X$ is a ring where the 
multiplication is defined by 
$$v\bullet u= \delta^! (v\times u)\,\,$$
$\delta^!$ being the refined Gysin map related to the diagonal 
embedding $X\to X\times X$. Hence we have: 
\smallskip 

\proclaim{5.8 Proposition} 
The multiplication of the ring $CH_{\dim X}X\times _S X$ 
respects $F^\bullet$. Namely if 
$u\in F^\nu CH_{\dim X}X\times_S X$ and $v\in F^\mu CH_{\dim X}
X\times _S X$, 
then $$v\bullet u\in F^{\nu+\mu}  CH_{\dim X}X\times _S X$$ 
In particular, if we assume the conjectures of Murre and of Grothendieck, 
the ideal 
\newline $F^1 CH_{\dim X}X\times _S X=CH_{\dim X}(X\times _S X)_{hom}$ 
is a nilpotent ideal.\qed\endproclaim 

\proclaim{5.9 Theorem} (Assume Grothendieck's and Murre's conjectures)
Let $f: X\to S$ be as above. Then the surjective map 
$$\rho: CH_{\dim X} X\times_S X \to 
\En_{\M S} {}^p \underline R f_* {Q}_X$$
has a nilpotent kernel.
\endproclaim

\demo{Proof}  
We claim first that the map of rings 
$$H^{BM}_{2\dim X} X\times_S X \to \En_{grPerv  S} 
\sum {}^p  R^i f_\ast {Q}_X[-i] $$
has a nilpotent kernel. If one chooses a decomposition
$$R f_* {Q}_{X} = \sum {}^p R^i f_* {Q}_{ X} [-i]$$
in $D^b_{cc} (S)$, an endomorphism 
$u$ of $R f_* {Q}_{X} $
may be represented by a matrix according to the decomposition. The matrix 
is upper triangular 
since the maps ${}^p R^i f_* {Q}_{ X} [-i] 
\to {}^p R^j f_* {Q}_{ X} [-j]$
are zero if $i<j$. 

Consider now an element 
$$u\in H^{BM}_{2\dim X} X\times_S X = \En_{D^b_{cc}  S} R f_* {Q}_X
$$  
which maps to zero in $ \En_{grPerv  S} \sum {}^pR^i f_\ast
{Q}_X[-i] $. 
Then the matrix representing $u$ as above is strictly upper triangular. 
Hence there exists $N$ such that $u^N=0$ as an endomorphism of  
$R f_* {Q}_{X}$. 

The kernel of the homomorphism 
$$CH_{\dim X} X\times _S X \to H^{BM} _{2\dim X} X\times_S X$$
equals $CH_{\dim X}(X\times _S X)_{hom}$.  Under the conjectures, 
it is a nilpotent ideal, 5.9.  Hence the kernel of $\rho $ is also 
nilpotent. \qed \enddemo 

\proclaim{5.10 Corollary} 
Any set of orthogonal projectors $\{\pi^i \}$ of $\En_{\M S} 
{}^p \underline Rf_* {Q}_X$ such that $\Delta_X = \sum \pi^i$ can be 
lifted to a set of orthogonal projectors $\{\Pi^i \}$ 
of \linebreak 
$CH_{\dim X} (X\times_S X)$ such that $\Delta_X = \sum \Pi^i$.\endproclaim

\demo{Proof} More generally the following holds (cf. [Ja1, Lemma 5.4]). \qed
\enddemo

\proclaim{5.11 Proposition}
Let  $\phi: A\to B$ is a surjective homomorphism of not necessarily 
commutative rings with nilpotent kernel.  Then any orthogonal set 
$\{ p_1, \cdots p_m\}$ of idempotents of $B$ (i.e. $p_i p_j = 
\delta_{i, j} p_i$) adding up to $1_B$ can be lifted to an orthogonal set 
of idempotents of $A$ adding up to $1_A$. \qed \endproclaim

\head 6. Decomposition in $CH\M (S)$ \endhead 

Throughout this section we assume that desingularisations of varieties over 
$k$ exist.
The aim is to prove the following:

\proclaim{6.1 Decomposition theorem in $CH \M (S)$}
Assume the standard conjectures and Murre's conjecture.
Let $X$ be a smooth
variety and $f: X \to S$ be a projective morphism.

(1) There is a noncanonical direct sum decomposition in $CH \M(S)$
$$\CR f_\ast {Q}_X \cong \sum_{m,k} 
\underline {\Cal C \IC} V^{m}_k[k-m]$$
where $V_k$ is a local system on a Zariski locally closed 
subvariety $T_k\subset S$ and 
$\underline {\Cal C \IC} V^{m}_k[k-m]$ denotes an object in $CH \M (S)$, 
together with a given 
isomorphism in $\De^b_{cc}(S)$
$$\r \underline {\Cal C \IC} V^{m}_k [k-m] 
@> \cong >> 
\IC V^{m}_k [k-m]$$

(2) The monomorphisms
$${}^p\tau_{\leq m} \CR f_\ast {Q}_X 
=\sum_{i\leq m}\sum_k \underline {\Cal C \IC} V^{i}_k[k-i]
\to \CR f_\ast {Q}_X$$
(where the 1st equality is the definition of ${}^p\tau_{\leq m} 
\CR f_\ast {Q}_X$)
are specified up to canonical isomorphism. In particular so are 
the ``subquotients''
$${}^p \CR^mf_\ast {Q}_X [-m]
=\sum_k \underline {\Cal C \IC} V^{m}_k[k-m]$$
(where the 1st equality is the definition of
${}^p \CR^mf_\ast {Q}_X [-m]$) specified up to canonical isomorphism.

(3) The decompositions
$${}^p \CR^mf_\ast {Q}_X [-m]= \sum_k \underline {\Cal C \IC} V^{m}_k[k-m]$$
are uniquely specified.
\endproclaim

\demo{Proof} Using 4.12, choose a decomposition in $\M (S)$
$${}^p \underline R f_\ast {Q}_X \cong \sum_{m,k} \underline \IC V^{m}_k[k-m]$$
and let 
$$\pi^m_k\in \En_{\M S}X$$ 
be the projector onto $\underline \IC V^{m}_k[k-m]$. By 5.10, 
the $\pi^m_k$ lift
to an orthogonal set of projectors 
$$\Pi^m_k \in \En_{CH \M S}X$$
Now set
$$\underline {\Cal C \IC} V^{m}_k[k-m] = (X, \Pi^m_k)$$
This proves the existence of the sought for decomposition. 

The uniqueness
statement (2) is more subtle, and will be deduced from the 
corresponding uniqueness 
statements for the decomposition in $\M (S)$. According to definition 
6.2 below,
${}^p\tau_{\leq m}\CR f_\ast {Q}_X$
(resp. ${}^p\tau_{>m}\CR f_\ast {Q}_X$) has cohomological
degree $\leq m$ (resp. $>m$). Then, by 6.3(2) below
$$\Ho_{CH \M S}\bigl({}^p\tau_{\leq m}\CR f_\ast {Q}_X,
{}^p\tau_{>m}\CR f_\ast {Q}_X\bigr)=0$$
{\it independently} of the choice of the liftings $\Pi^m_k$. This implies (2).

Finally, to prove the uniqueness statement in (3), note that by construction
${}^p \CR^mf_\ast {Q}_X [-m]$ has cohomological degree exactly $m$, hence
the statement follows from the decomposition theorem in $\M (S)$ 4.12, and
6.3(2). \qed \enddemo

The rest of this section is devoted to the proof of 6.3, which was 
used in the proof of 6.1.

\definition{6.2 Definition} A Chow motive $(X, P)$ over $S$ has cohomological
degree $\leq m$ (resp. $>m$) if
$${}^p {\Cal H}^i \r (X,P)=0$$
for all $i>m$ (resp. $i\leq m$).
Finally, $(X, P)$ has degree exactly $m$ if it has degree $\geq m$ and 
$\leq m$.
\enddefinition

\proclaim{6.3 Theorem} (Assuming the standard and Murre's conjectures)

(1) If $(X,P)$ has cohomological degree $\leq m$ and
$(Y,Q)$ has cohomological degree $>m$, then 
$$\Ho_{CH \M S}\bigl( (X,P), (Y,Q) \bigr)=0$$

(2) If $(X,P)$ and $(Y,Q)$ have degrees exactly $m$, then
$$\Ho_{CH \M S}\bigl( (X,P), (Y,Q) \bigr)=\Ho_{\M S}
\bigl( (X,P), (Y,Q) \bigr)
$$
\endproclaim

\demo{Proof} Let $Z=X\times_S Y$. We consider the operators
$$CH_s Z \ni \alpha \to \Psi \alpha = Q \bullet \alpha \bullet P \in CH_s Z$$
$$H^{BM}_i Z \ni a \to \psi a =[Q]\bullet a \bullet [P] \in H^{BM}_i Z$$
Note that
$\Psi^2=\Psi$ and $\psi^2 =\psi$, i.e. both operators are projectors. 

To prove (1) and (2), it is enough to show that if
$\psi =0$, then $\Psi=0$. Most of the argument
will be spent showing that, because $\psi =0$, $\Psi F^\nu \subset F^{\nu+1}$
for all $\nu$.

{\smc Step 1.} Let $T \subset Z$ be the singular set. Make a diagram
$$\diagram
\overline R \rto & \overline W\\
R \uto\rto \dto  & W \uto \dto^p\\
T\rto            & Z \enddiagram$$
where: 

(a) $\overline W$ is smooth and projective and $\overline R \subset \overline W$
is a smooth normal crossing divisor,

(b) $(W, R)=(\overline W -B, \overline R -B)$ for some divisor $B \subset 
\overline R$, 

(c) $p: W \to Z$ is a resolution of singularities and $R=p^{-1} T \subset W$.

By 6.8 below, $\Psi$ is a class $\C$ operator, hence by 6.5 there is
a correspondence
$$\Gamma \in CH_\bullet (\overline W \times \overline W)$$
($\overline W$ is not necessarily equidimensional)
such that
$$\Psi \bigl( p_\ast (\alpha |W)\bigr)=p_\ast\bigl((\Gamma_\ast \alpha )|W
\bigr)\eqno (\ast)$$
for all $\alpha \in CH_\bullet Z$, resp.
$$\psi \bigl( p_\ast (a |W)\bigr)=p_\ast\bigl((\Gamma_\ast a)|W
\bigr)\eqno (\ast \ast)$$
for all $a\in H^{BM}_\bullet Z$.

{\smc Step 2.} We have an exact sequence
$$H^{BM}_i \overline R \to H^{BM}_i R \oplus H^{BM}_i
\overline W  \to H^{BM}_i W$$

If $j : W \hookrightarrow \overline W$, $j^\prime : R \hookrightarrow \overline
R$ are the natural inclusions, the exact sequence arises, by means of a 
familiar construction, from the following morphism of distinguished
triangles
$$\diagram
D_B \ddouble \rto & D_{\overline R}\dto \rto &Rj^\prime_\ast D_R \dto
\rto^{[1]}&\\
D_B \rto          & D_{\overline W}\rto  &Rj_\ast D_W\rto^{[1]}&\enddiagram$$

{\smc Step 3.} This is the crucial step. We show that
$$\gamma =cl \Gamma \in N^1 H_\bullet^{BM} \overline W \times \overline W$$
lies in the codimension 1 piece of the coniveau filtration.

Let $i: \overline R \hookrightarrow \overline W$ and
$i^\prime : R \hookrightarrow W$ be the inclusions.
First of all, from the standard exact sequence
$$H_i^{BM} R \to H_i^{BM} T \oplus H_i^{BM} W \to H_i^{BM} Z$$
we deduce that if
$$p_\ast (a|W)=0$$
for some $a \in H_i^{BM} \overline W$, then
$$a|W= i_\ast^\prime a^\prime$$
for some $a^\prime \in H_i^{BM} R$. Then, from the exact sequence in 
Step 2, 
$$a =i_\ast a^{\prime \prime}$$
for some $a^{\prime \prime} \in H^{BM}_i \overline R$. 
Now, for any $a \in H^{BM}_i \overline W$
$$p_\ast \bigl(( \Gamma_\ast a)|W\bigr)
=\psi \bigl( p_\ast (a |W)\bigr)=0$$
by assumption, so, because of what has been said
$$\Gamma_\ast a \in \Img H^{BM}_\bullet \overline R$$
which finishes step 3.

{\smc Step 4.} Let $\Gamma^\prime = \Gamma \circ \Gamma \circ ...$
(many times). Since $\Psi^2=\Psi$ and $\psi^2 =\psi$, equations
$(\ast)$ and $(\ast \ast)$ are still satisfied with $\Gamma^\prime$ in
place of $\Gamma$. In addition we have
$$\gamma^\prime =cl \Gamma^\prime \in N^{\text{many times}}
H_\bullet \overline W \times \overline W =0$$
Then $\Gamma_\ast^\prime F^\nu CH^\bullet \overline W \subset 
F^{\nu + 1} CH^\bullet \overline W$ for all $\nu$ and, by the 
strictness properties of the $F$-filtration
$$\Psi F^\nu CH_\bullet Z \subset F^{\nu + 1} CH_\bullet Z$$
for all $\nu$. 

{\smc Step 5.} We are assuming the standard conjectures. Therefore
$F^\nu CH_\bullet Z=0$ for $\nu$ large. This implies that $\Psi =0$,
which concludes the proof.\qed \enddemo

The rest of the paper is devote to finishing the proof of 6.3.

\definition{6.4 Definition}

(1) Let $X$ be a quasiprojective variety. For the purpuse of the following
discussion, a {\it smooth cover} of $X$ is a diagram
$$\diagram
U \dto^p \rto& \overline U \\
X & \enddiagram$$
sometimes simply denoted $\overline U \supset U \overset p \to \to X$,
where $\overline U$ is a smooth projective variety, $U \subset
\overline U$ an open subvariety with smooth normal crossing boundary
divisor $\overline U \setminus U$, and $p : U \to X$ a projective 
morphism. It is a consequence of our assumptions on existence of 
resolution of singularities, that smooth covers always exist.

(2) Let $X$, $Y$ be quasiprojective varieties. We say that an operator
$$\Psi : CH_\bullet X \to CH_{\bullet -c} Y$$
is of {\it class} $\C$ ($\C$ stands for ``correspondence'') if there 
are smooth covers $\overline U \supset U \overset p \to \to X$ 
of $X$ and $\overline V \supset V \overset q \to \to Y$ of $Y$, and
a correspondence
$\Gamma \in CH_\bullet \overline U \times \overline V$
such that
$$\Psi \bigl( p_\ast(\alpha |U)\bigr)
= q_\ast \bigl( (\Gamma_\ast \alpha ) |V\bigr)$$
for all $\alpha \in CH_\bullet \overline U$. In this case, we say that 
$\Gamma$ {\it induces} $\Psi$.
\enddefinition

The following is the basic point

\proclaim{6.5 Lemma} Let $\Psi : CH_\bullet X \to 
CH_{\bullet -c} Y$ be of class $\C$, and
let $\overline U \supset U \overset p \to \to X$, 
$\overline V \supset V \overset q \to \to Y$ be arbitrary smooth covers of
$X$, $Y$. Then, there exists a correspondence
$\Gamma \in CH_\bullet \overline U \times \overline V$ inducing $\Psi$
\endproclaim

\demo{Proof} By assumption, there are {\it some} smooth covers
$\overline U_1 \supset U_1 \overset p_1 \to \to X$ of $X$, 
$\overline V_1 \supset V_1 \overset q_1 \to \to Y$ of $Y$,
and {\it some} correspondence
$\Gamma_1 \in CH_\bullet \overline U_1 \times \overline V_1$
such that
$$\Psi \bigl( p_{1\ast}(\alpha |U_1)\bigr)
= q_{1\ast} \bigl( (\Gamma_{1\ast} \alpha )|V_1)\bigr)$$
for all $\alpha \in CH_\bullet \overline U_1$.

Now, any 2 smooth covers can be housed under a third
$$\diagram
    &U_2 \subset \overline U_2 \dlto \drto&      \\
U_1 \subset \overline U_1& &U \subset \overline U\enddiagram$$
so, in the end, we may assume that there is either 
a morphism $U_1 \subset \overline U_1 \to U \subset \overline U$,
or the other way around, and similarly for $V$.

{\smc Case 1.} Assume that there is a morphism $\pi : 
U \subset \overline U \to U_1 \subset \overline U_1$ and let
$$\Gamma  =\Gamma_1 \circ \Gamma_\pi \in CH_\bullet 
(\overline U\times \overline V_1)$$
Then $\Gamma$ induces $\Psi$, since
$$\align
\Psi \bigl(p_\ast (\alpha|U)\bigr)&=\Psi\bigl(p_{1\ast}\pi_\ast (\alpha|U)
\bigr)=\Psi\bigl(p_{1\ast} ((\pi_\ast \alpha)|U_1)\bigr)=\\
&=q_{1\ast}\bigl((\Gamma_{1\ast} \pi_\ast \alpha)|V_1\bigr)=
q_{1\ast}\bigl( (\Gamma_\ast \alpha)|V_1\bigr)\endalign$$

{\smc Case 2.} Assume now that there is a morphism $\pi : 
U_1 \subset \overline U_1 \to U \subset \overline U$. Let 
$i : \overline W \hookrightarrow \overline U_1$ be a smooth projective
subvariety generically finite over $\overline U$
$$\diagram
\overline W \rto^i & \overline U_1 \dto^\pi \\
                   & \overline U \enddiagram$$
and let us agree that $d$ be the generic degree of $\pi \circ i$.
Let us fix ourselves a correspondence $\Gamma \in CH_\bullet \overline U 
\times \overline V_1$ with the property that
$$\Gamma _ \ast = {1\over d} \Gamma_{1 \ast}i_\ast i^\ast \pi^\ast$$
Then $\Gamma$ induces $\Psi$, since
$$\align
\Psi\bigl( p_\ast (\alpha|U)\bigr)&
=\Psi\bigl(p_\ast\bigl({1\over d} \pi_\ast i_\ast i^\ast \pi^\ast \alpha|
U\bigr)\bigr)=\\
&=\Psi\bigl(p_{1\ast}\bigl({1\over d} i_\ast i^\ast \pi^\ast \alpha|
U\bigr)\bigr)=\\
&=q_{1\ast}\bigl(\bigl({1\over d}\Gamma_{1\ast}i_\ast i^\ast\pi^\ast
\alpha\bigr)|V_1\bigr)=q_{1\ast}\bigl((\Gamma_\ast \alpha)|V_1\bigr)\endalign$$

To summarise, in both cases we were able to find a $\Gamma\in CH_\bullet
\overline U \times \overline V_1$ inducing $\Psi$. Working now
the $V$s in a similar fashion, we can also find
$\Gamma \in \overline U \times \overline V$ inducing $\Psi$, i.e., prove
the lemma. \qed \enddemo

\proclaim{6.6 Lemma} Let $\Psi : CH_\bullet X \to CH_{\bullet -c} Y$ and
$\Phi : CH_\bullet Y \to CH_{\bullet -e} Z$ be of class $\C$. 
Then, the composition
$\Phi \circ \Psi: CH_\bullet X \to CH_{\bullet -c-e} Z$ is also of class $\C$.
\endproclaim

\demo{Proof} Obvious. \qed \enddemo

\proclaim{6.7 Lemma}

(1) Let $f:X \to Y$ be a proper map. Then $f_\ast :CH_iX\to CH_iY$ is
of class $\C$.

(2) Let $i:Y \hookrightarrow X$ be a regular embedding of codimension
$c$ and 
$$\diagram 
Y^\prime \dto \rto & X^\prime \dto \\
Y \rto^i           &X \enddiagram$$
be a fibre square. The refined Gysin map
$$i^!:CH_\bullet X^\prime \to CH_{\bullet-e} Y^\prime$$
is of class $\C$
\endproclaim

\demo{Proof} (1) is obvious and (2) is 5.6. \qed \enddemo

\proclaim{6.8 Corollary} The operator $\Psi$ in the proof of 6.3 is
of class $\C$.
\endproclaim

\demo{Proof} Indeed, $\Psi$ is a composition of proper push forward and
refined Gysin maps. All these are of class $\C$ 6.7, and so their
composition is 6.6.\qed \enddemo

\enddocument